\documentclass[processes,article,submit,moreauthors,pdftex,12pt,a4paper]{mdpi} 
\pdfoutput=1
\usepackage{subfigure}
\lastpage{x}
\doinum{10.3390/------}
\pubvolume{xx}
\pubyear{2013}
\history{Received: xx / Accepted: xx / Published: xx}

\Title{Reduced models in chemical kinetics via nonlinear data-mining}

\Author{Eliodoro Chiavazzo $^{1,2}$, C. William Gear $^{1}$, Carmeline J. Dsilva $^{1}$, Neta Rabin $^{3}$ and Ioannis G. Kevrekidis $^{1,4,}$*}

\address{%
$^{1}$ Department of Chemical and Biological Engineering, Princeton University, Princeton, NJ USA\\
$^{2}$ Energy Department, Politecnico di Torino, Torino, Italy\\
$^{3}$ Department of Exact Sciences, Afeka Tel-Aviv Academic College of Engineering, Tel-Aviv, Israel\\
$^{4}$ Program in Applied and Computational Mathematics, Princeton University, Princeton, NJ USA
}

\corres{yannis@princeton.edu; Tel: 609 258-2818; Fax: 609 258-0211}

\abstract{The adoption of detailed mechanisms for chemical kinetics often poses two types of severe challenges:
First, the number of degrees of freedom is large; and second, the dynamics is characterized by widely disparate time scales.
As a result, reactive flow solvers with detailed chemistry often become intractable even for large clusters of CPUs, especially when dealing with direct numerical simulation (DNS) of turbulent combustion problems.
This has motivated the development of several techniques for reducing the complexity of such kinetics models, where eventually only a few variables are considered in the development of the simplified model.
Unfortunately, no generally applicable {\em a priori} recipe for selecting suitable parameterizations of the reduced model is available, and the choice of slow variables often relies upon intuition and experience.
We present an automated approach to this task, consisting of three main steps. First, the low dimensional manifold of slow motions is (approximately) sampled by brief simulations
of the detailed model, starting from a rich enough ensemble of admissible initial conditions.
Second, a {\em global} parametrization of the manifold is obtained through the Diffusion Map (DMAP) approach, which has recently emerged as a powerful tool in data analysis/machine learning.
Finally, a simplified model is constructed and solved {\em on the fly} in terms of the above reduced (slow) variables.
Clearly, closing this latter model requires nontrivial interpolation
calculations, enabling restriction (mapping from the full ambient space to the reduced one) and lifting (mapping from the reduced space to the ambient one).
This is a key step in our approach, and a variety of interpolation schemes are reported and compared.
The scope of the proposed procedure is presented and discussed by means of an illustrative combustion example.}


\begin{document}

\section{Introduction}
The solution of detailed models for chemical kinetics is often dramatically time consuming owing to a large number of variables evolving in processes with a wide range of time and space scales. As a result, fluid dynamic flow solvers coupled with detailed chemistry still present a challenge, even for modern clusters of CPUs, especially when dealing with direct numerical simulation (DNS) of turbulent combustion systems. Here, a large number of governing equations for chemical species (a few hundred for mechanisms of standard hydrocarbon fuels) are to be solved at (typically) millions of distinct discretization points in the computational domain. This has motivated the development of a plethora of approaches aiming at reducing the computational complexity of such detailed combustion models, ideally by recasting them in terms of only a few new reduced variables. (see e.g. \cite{MaasGoussis} and references therein).
The implementation of many of these techniques typically involves three successive steps. First, a large set of stiff ordinary differential equations (ODEs) is considered for modeling the temporal evolution of a spatially homogenous mixture of chemical species under specified stoichiometric and thermodynamic conditions (usually fixed total enthalpy and pressure for combustion in the low Mach regime). It is well known that, due to the presence of fast and slow dynamics, the above systems are characterized by {\em low dimensional manifolds} in the concentration space (or phase-space), where a typical solution trajectory is initially rapidly attracted towards the manifold, while afterwards it proceeds to the thermodynamic equilibrium point always remaining in close proximity to the manifold. Clearly, the presence of a manifold forces the ODEs state to visit mostly a low dimensional region of the entire phase-space, thus offering the premise for constructing a consistent reduced description of the process, which accurately retains the slow dynamics along the manifold while neglecting the initial short transient towards the manifold.
In a fluid dynamic simulation, stoichiometry and thermodynamic conditions may vary throughout the computational domain. Hence, when implementing reduction techniques, the second step consists of parameterizing and tabulating the manifolds arising in the homogeneous reactor for a variety of stoichiometric and thermodynamic conditions. Finally, as a third step, the fluid dynamic equations are reformulated in terms of the new variables, with the latter tables utilized to close
the new reduced set of equations (see, e.g., \cite{ChiavazzoCF}).
It is worth stressing that the above description briefly outlines only one possible approach for coupling a model reduction method to a flow solver: the case
 where the low dimensional manifolds of the homogeneous problem are identified in advance {\em in the entire phase-space}.
For completeness, it is important mentioning that, due to the rapidly increasing difficulty in storing and interpolating data in high dimensions, this approach remains viable in cases with {\em a few} reduced variables. As an alternative to this global method, techniques have been introduced for {\em locally} constructing the low dimensional manifold only in the (tiny) region of interest in the phase-space, as demanded by a reacting flow code during simulations \cite{ICEPIC,ChiavazzoKarlinPRE,Chiavazzo2012}.
Local constructions can certainly cope with higher dimensional manifolds. However, their usage seems computationally advantageous only in combination with efficient algorithms for {adaptive} tabulation, where data is computed when needed, stored, and re-utilized if necessary (see, e.g., \cite{ISAT}). 

%
In this work, we focus on the global construction and parameterization of slow invariant manifolds arising in the modeling of spatially homogeneous reactive mixtures. In particular, upon identification of the slow manifold, we propose a generally applicable methodology for selecting a suitable parameterization; we also investigate various interpolation/extrapolation schemes that need to be used in the solution of a reduced dynamical system expressed in terms of the variables learned.

The manuscript is organized as follows. In Section \ref{DMAP}, Diffusion Maps are briefly reviewed. In Section \ref{approach} and subsections therein, we discuss the computation of points on the manifold, their embedding in a reduced (here two-dimensional) space, the formulation of a reduced set of equations and their solution through several interpolation/extension techniques.
Results are reported and discussed in Section \ref{results}, where the proposed approach is applied to a reactive mixture of hydrogen and air at stoichiometric proportions with fixed enthalpy and pressure. The reader may prefer a quick glance at Section \ref{results} before the detailed presentation of the procedure in Section \ref{approach}.
Finally, we conclude with a summary and brief discussion of open issues in Section \ref{conclusion}.

\section{Diffusion maps}\label{DMAP}
The Diffusion Map (DMAP) approach has emerged as a powerful tool in data analysis and dimension reduction \cite{Coifman05pnas01,Coifman05pnas02,Coifman2006}. In effect, it can be thought of as a nonlinear counterpart of Principal Component Analysis (PCA) \cite{PCAref} that can be used to search for a low-dimensional embedding of a high-dimensional point set \{${\bf y}_1$,...,${\bf y}_M$\}, if an embedding exists. For completeness, we present a simple description of the DMAP process.
The points ${\bf y}_i$ could exist in some $n$-dimensional Cartesian space (as they are in our combustion example) or they could be more abstract objects, such as images. What is important is that there exists a {\em dissimilarity function}, $d_{ij} = d_{ji}$ between any pair of points, ${\bf y}_i$ and ${\bf y}_j$ such that the dissimilarity is zero only if the points are identical (in those aspects that are important to the study) and gets larger the more dissimilar they are. Although, for points in $\Re^n$, an obvious choice for $d_{ij}$ is the standard Euclidean distance, this is not necessarily the best option. For instance, a weighted Euclidean norm may be considered when different coordinates are characterized by disparate orders of magnitude. As discussed below, this is indeed the case encountered in many combustion problems, where the data are composition vectors in concentration space and major species (i.e. reactants and products) are characterized by much higher concentrations compared to minor species (i.e. radicals).  From $d_{ij}$ a pairwise affinity function $w_{ij}=w(d_{ij})$ is computed where $w(0) = 1$ and $w(d)$ is monotonically decreasing and non-negative for $d > 0$. A popular option is the heat kernel
\begin{equation}\label{similarity}
w \left( d \right) = \exp \left[ { - \left( {\frac{{d }}{\varepsilon }} \right)^2 } \right].
\end{equation}
The model parameter $\varepsilon$ specifies the level below which points are considered similar, whereas points more distant than a small multiple of $\varepsilon$ are, effectively, not linked directly. For this presentation we will assume that $d$ is a distance measure in (suitably scaled) Cartesian coordinates so that each point, ${\bf y}_i$ is specified by its coordinates, $y_{i,\alpha}$ with $\alpha=1,...,n$ in $n$-dimensional space.

%
In the DMAP approach, starting from the $M \times M$ (not $n \times n$ as in PCA) symmetric matrix $W=\{w_{ij}\}$, a Markov matrix $K$ is constructed through the row normalization
\begin{equation}\label{Markov}
K=D^{-1} W,
\end{equation}
with the diagonal matrix $D$ collecting all the row sums of matrix $W$. Owing to similarity with a symmetric matrix, $D^{-1/2} W D^{-1/2}$, $K$ has a complete set of real eigenvectors $\{ \phi_i \}$ and eigenvalues $\{ \lambda_i \}$.
%
%
%
%
%
%
%
Moreover, a projection of the high-dimensional points \{${\bf y}_1$,...,${\bf y}_M$\} into an $m$-dimensional space (hopefully $m << n$) can be established through the components of $m$ 
appropriately selected eigenvectors (not necessarily the $m$ leading ones, as in PCA).
Specifically, let the eigenvalues be sorted in decreasing order: $1=\lambda_1 \ge \left| {\lambda _2 } \right| \ge ... \ge \left| {\lambda _M} \right|$. 
The diffusion map $\Psi_t$ is defined based on the right eigenvectors of $K$, $K  \phi_l = \lambda_l  \phi_l$, with $\phi _l  = \left( {\phi _{1,l} ,...,\phi _{M,l} } \right)$, for $t>0$, as follows:
\begin{equation}\label{DMAPdef}
{\rm \Psi} _t \left( { {\bf y}_i } \right) = \left( {\begin{array}{*{20}c}
   {\lambda _1^t  \phi _{i,1}}  \\
   {\lambda _2^t  \phi _{i,2}}  \\
    \vdots   \\
   {\lambda _M^t \phi _{ i,M }}  \\
\end{array}} \right),
\end{equation}
and it assigns a vector of $M$ new coordinates to each data point ${\bf y}_i$.
Notice that all points have the same first coordinate in (\ref{DMAPdef}), since $ \phi_1$ is proportional to the all-ones vector (with eigenvalue 1).
Notice that the diffusion map coordinates are time-dependent; using longer times in the diffusion process damps high frequency components, so that
fewer coordinates suffice for an approximation of a given accuracy.
%
%
However, in order to achieve a drastic dimension reduction, for a fixed threshold $0<\delta<1$, it is convenient to define a {\em truncated} diffusion map:
\begin{equation}\label{DMAPdefTR}
\Psi _t^\delta  \left( {{\bf y}_i } \right) = \left( {\begin{array}{*{20}c}
   {\lambda _2^t  \phi _{i,2 }}  \\
   {\lambda _3^t  \phi _{i,3}}  \\
    \vdots   \\
   {\lambda _{m+1}^t  \phi _{ i,{m+1}  }}  \\
\end{array}} \right)
\end{equation}
where $m+1$ is the largest integer for which $\left| {\lambda _{m+1} } \right|^t  > \delta $.
Below we will consider only the eigenvector entries (i.e. take $t=0$), and will separately discuss using the eigenvalues (and their powers)
to ignore noise.

If the initial data points \{${\bf y}_1$,...,${\bf y}_M$\} are located on a (possibly non-linear) low dimensional manifold with dimension $m$, one might expect (by analogy to PCA) 
that a procedure exists to systematically select $m$ diffusion map eigenvectors for embedding the data.
%
%
%
%
%

If the points are fairly evenly distributed across the low-dimensional manifold, it is known that the principal directions of the manifold are spanned by some of the leading
eigenvectors (i.e., those corresponding to larger eigenvalues) of the DMAP operator and the corresponding eigenvalues are approximately
\begin{equation}\label{eigenvalK}
\lambda = 1 - \delta [k \pi d/L_\alpha]^2        
\end{equation}
where $\delta \approx \exp(-{d/ \varepsilon}^2)$, $d$ is the typical spacing between neighbors, and $L_\alpha$ is the length of the $\alpha$-th principal direction. 
Here $k = 1, 2, \cdots$ indicates the successive harmonics of the eigenvectors. 
(This approximation can be obtained by considering the regularly-spaced data case, assuming that $\varepsilon$ is comparable to $d$, and that $\delta$ is small enough that higher powers can be ignored.)
Section \ref{DMAPissues} below discusses how to ignore eigenvectors that are harmonics of previous ones by checking for dependence.  
Eq. (\ref{eigenvalK}) provides a tool for deciding when to ignore the smaller eigenvalues.  
Suppose, for example, that we know that our data accuracy is approximately a fraction $\gamma$ of the range of the data. 
This range roughly corresponds to the longest principal direction, say $L_1$.   
There is little point in considering manifold directions of the order of $\gamma L_1$, since they are of the order of the errors in the data.  
Hence by applying (\ref{eigenvalK}) we should ignore any eigendirections whose eigenvalue is less than $1 - (1 - \lambda_2)\gamma^{-2}$,
where $\lambda_2$ is the first non-trivial eigenvalue.

%
\subsection{Issues in the implementation of the algorithm.}\label{DMAPissues}
\begin{figure}[ht]
 \centering
 \subfigure{
  \includegraphics[scale=0.25]{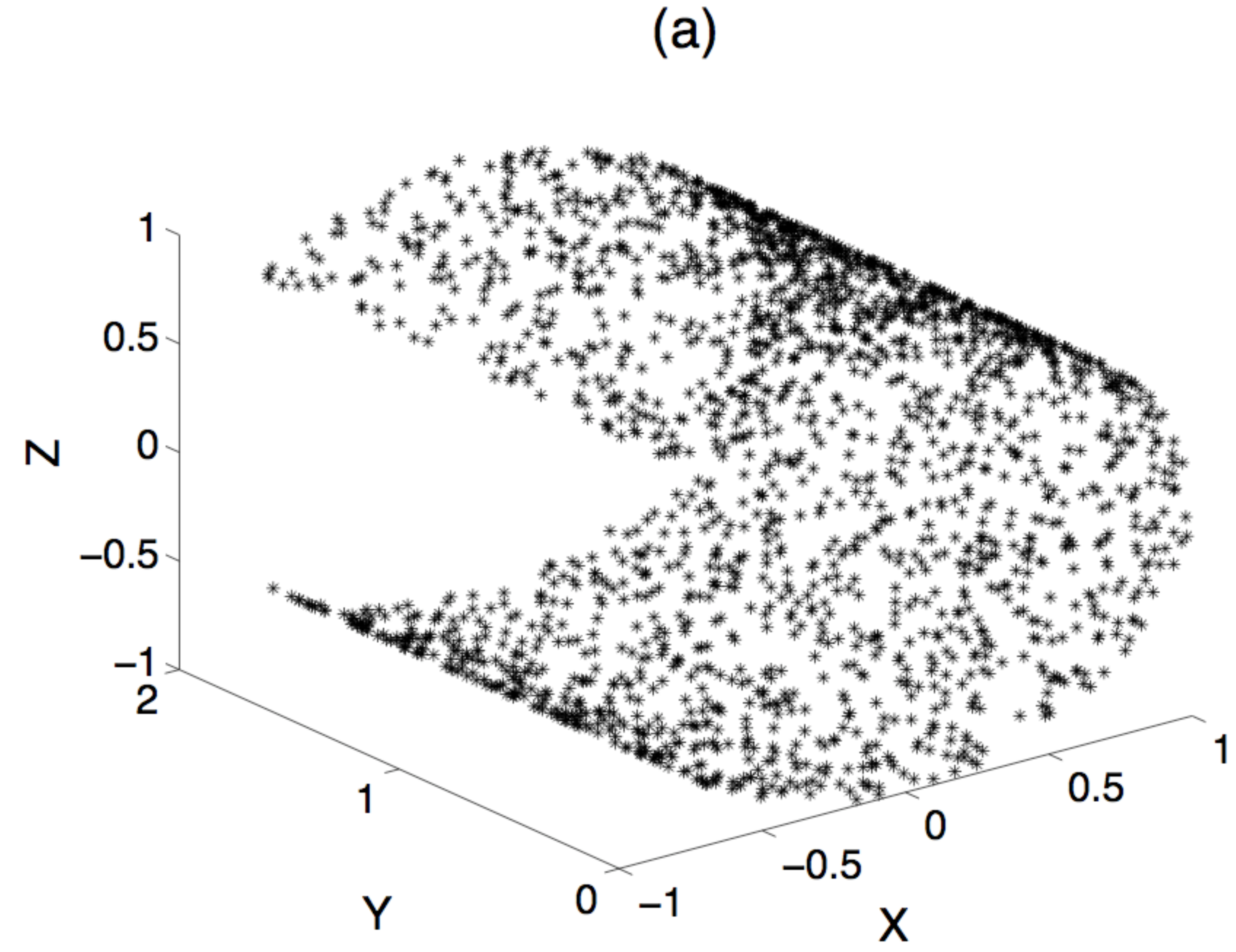}
   }
 \subfigure{
  \includegraphics[scale=0.25]{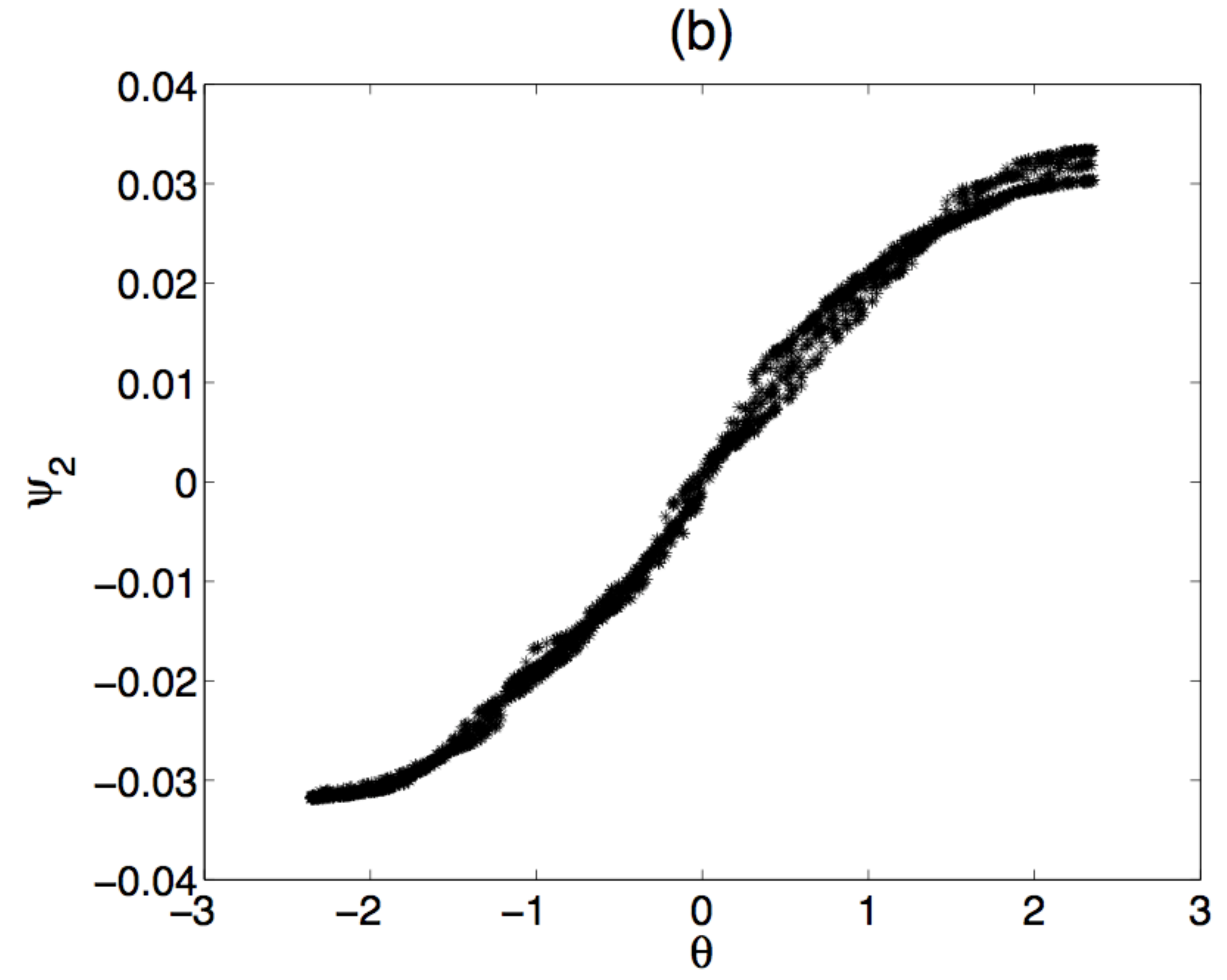}
   }
    \subfigure{
  \includegraphics[scale=0.25]{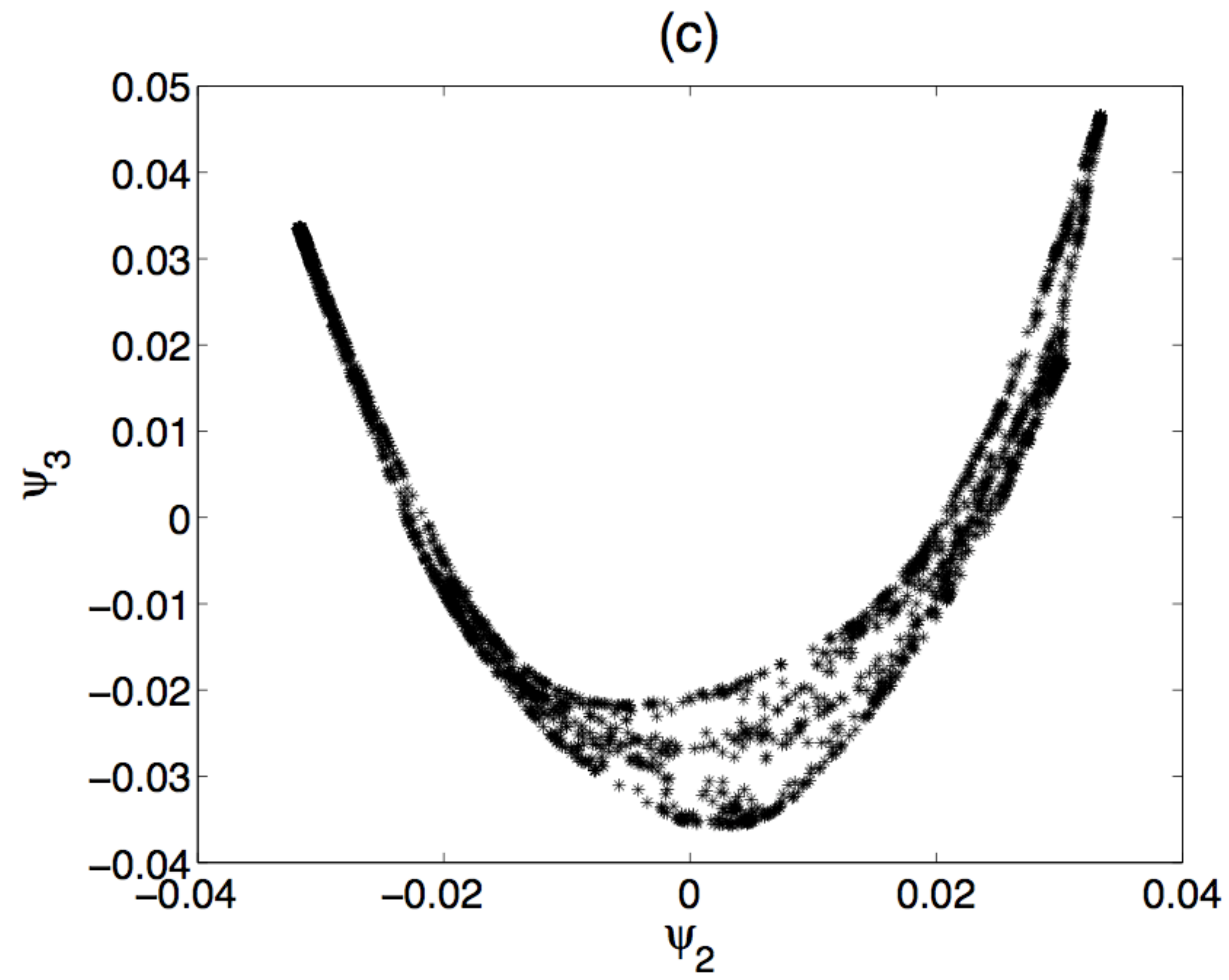}
   }
    \subfigure{
  \includegraphics[scale=0.25]{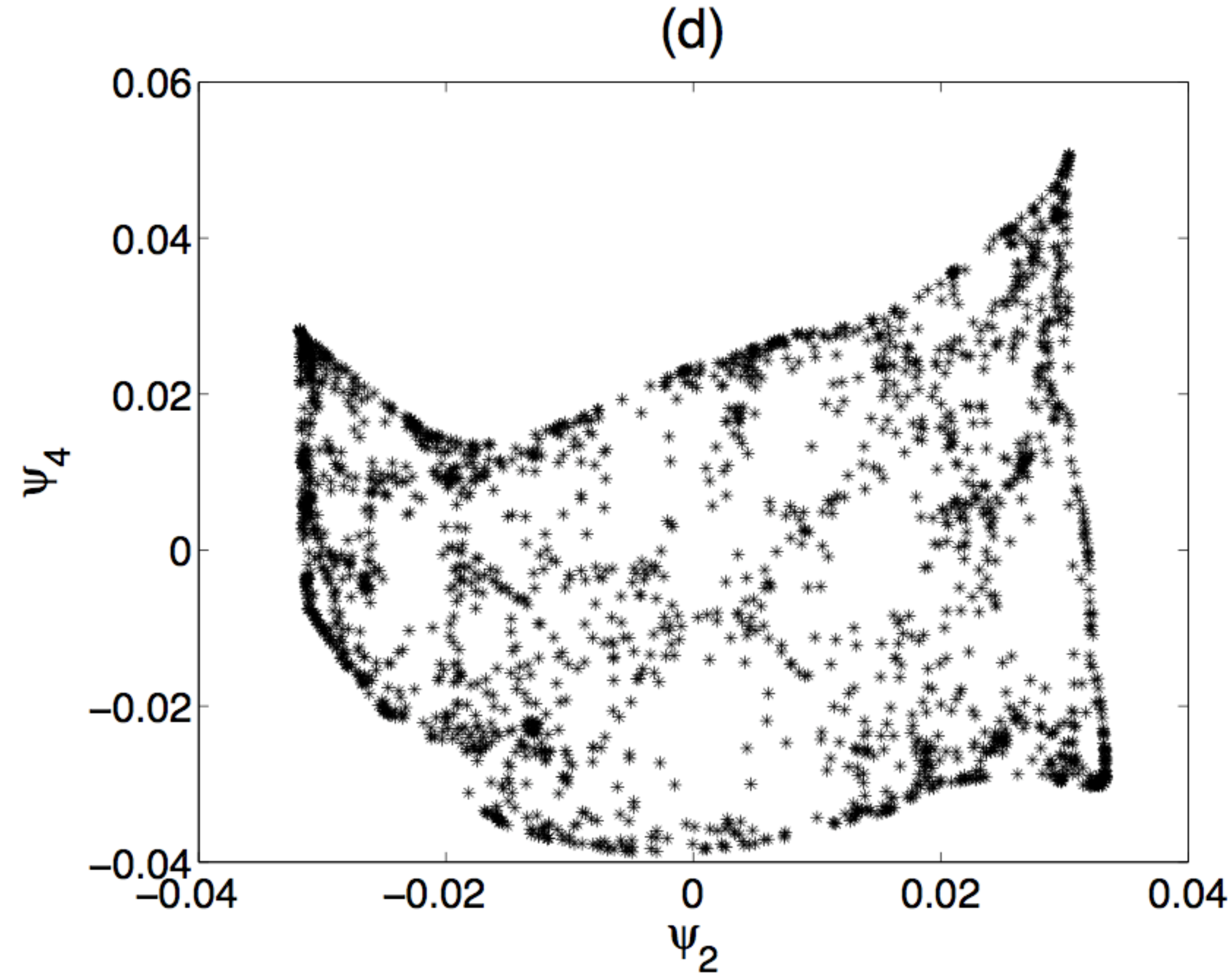}
   }
 \caption[Optional caption for list of figures]{%
Data manifold, dimensionality and independence of DMAP eigenvectors (a) 2000 uniformly random points initially placed in a unit square are stretched and wrapped around three fourths of a cylinder. (b) The entry in the first non-trivial eigenvector of the Markov matrix $K$ versus the first cylindrical coordinate $\theta$ for each
data point. (c) Entry in the second non-trivial eigenvector of $K$ versus the first one; the quasi-one-dimensionality of the plot implies strong eigenvector correlation. (d) Entry in the third non-trivial eigenvector of $K$ versus the first one. The evident two-dimensional scatter implies that a new direction on the data manifold has been detected.}\label{CYL01}
\end{figure}

While the formulas above appear to provide a simple recipe, a number of important, problem-dependent issues  arise,
having to do with the sampling of the points to be analyzed, the choice of the parameter $\varepsilon$ etc.; we
now discuss these issues through illustrative caricatures.
Consider 2000 uniformly random points initially placed in a unit square, then stretched and wrapped around three fourths of a cylinder of radius 1 and length 2 (see Fig. \ref{CYL01}(a)).
In Fig. \ref{CYL01}(b) the first non-trivial eigenvector, $ \psi_2$, is reported against the first cylindrical coordinate $\theta$:
the i-th component of this eigenvector is plotted against the $\theta$ angle of the i-th point.
The clearly apparent one-dimensional nature of the plot confirms that $ \psi_2$ parametrizes this principal geometric direction.
However, a plot of the $ \psi_3$, the eigenvector corresponding to the next leading eigenvalue, against $ \psi_2$ clearly shows a strong correlation:
$ \psi_3$ is not representative of a new, {\em indepedent} direction on the data manifold.
In Fig. \ref{CYL01}(d), the two-dimensional scatter of the plot of the entries of the fourth eigenvector versus the entries of the second one
indicates independence between $ \psi_2$ and $ \psi_4$; $ \psi_4$ does represent a new, independent direction along the data manifold and
becomes our second embedding coordinate.
Visually testing independence between two DMAP eigenvectors is relatively easy: we can agree that Figs. \ref{CYL01}(b)(c) appear one-dimensional and
Fig. \ref{CYL01}(d) appears two-dimensional.
But testing independence in higher dimensions (for subsequent DMAP eigenvectors) becomes quickly visually impossible and even computationally nontrivial.
Subsequent eigenvectors should be plotted against $\psi_2$ and $ \psi_4$ and the dimensionality of the plot should be assessed;
this is still visually doable for, say, $ \psi_5$, and the plot appears as a 2-D surface in 3-D:  $ \psi_5$ is not a new data coordinate.
Beyond visual assessment (and in higher dimensions) one can use the sorted edge-length algorithm for dimensionality assessment:
a log-log plot of the graph edge-length versus edge number is constructed, with the manifold dimension being the slope in the middle part of the plot.
Algorithms for detecting the dimension of attractors in chaotic dynamical systems can also find use here \cite{Grass81,GrassProc83}.
%

\begin{figure}[ht]
 \centering
 \subfigure{
  \includegraphics[scale=0.25]{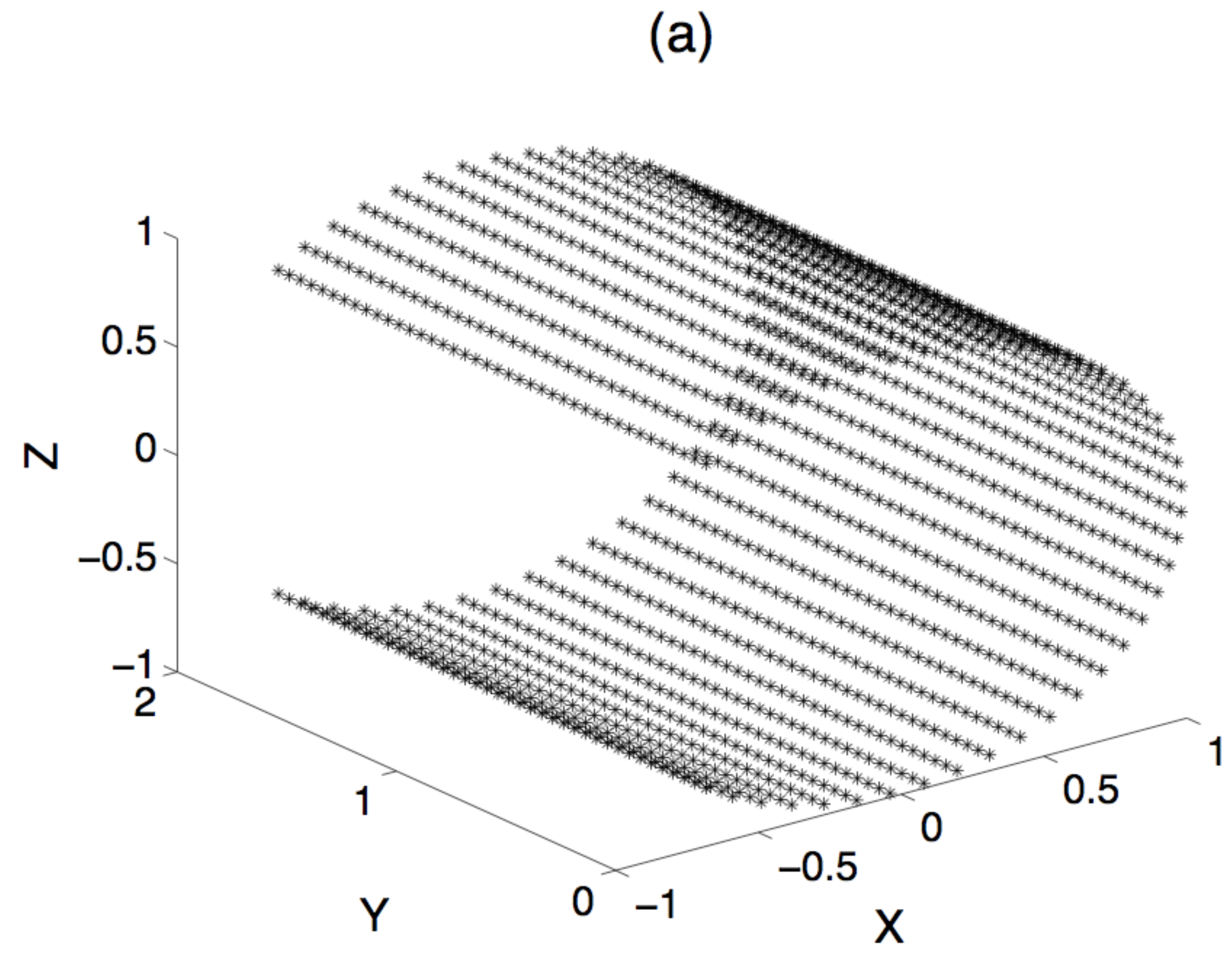}
   }
 \subfigure{
  \includegraphics[scale=0.25]{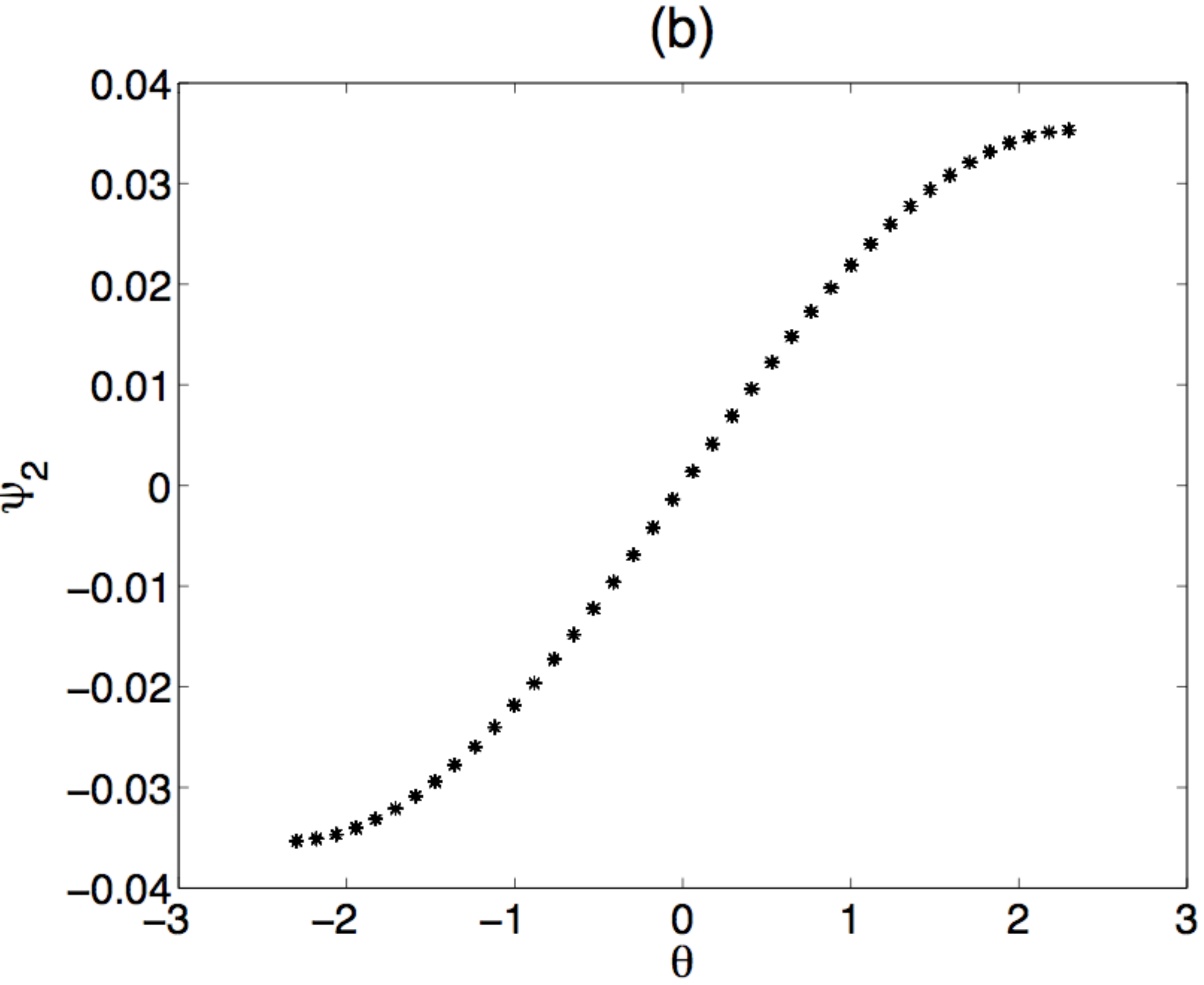}
   }
    \subfigure{
  \includegraphics[scale=0.25]{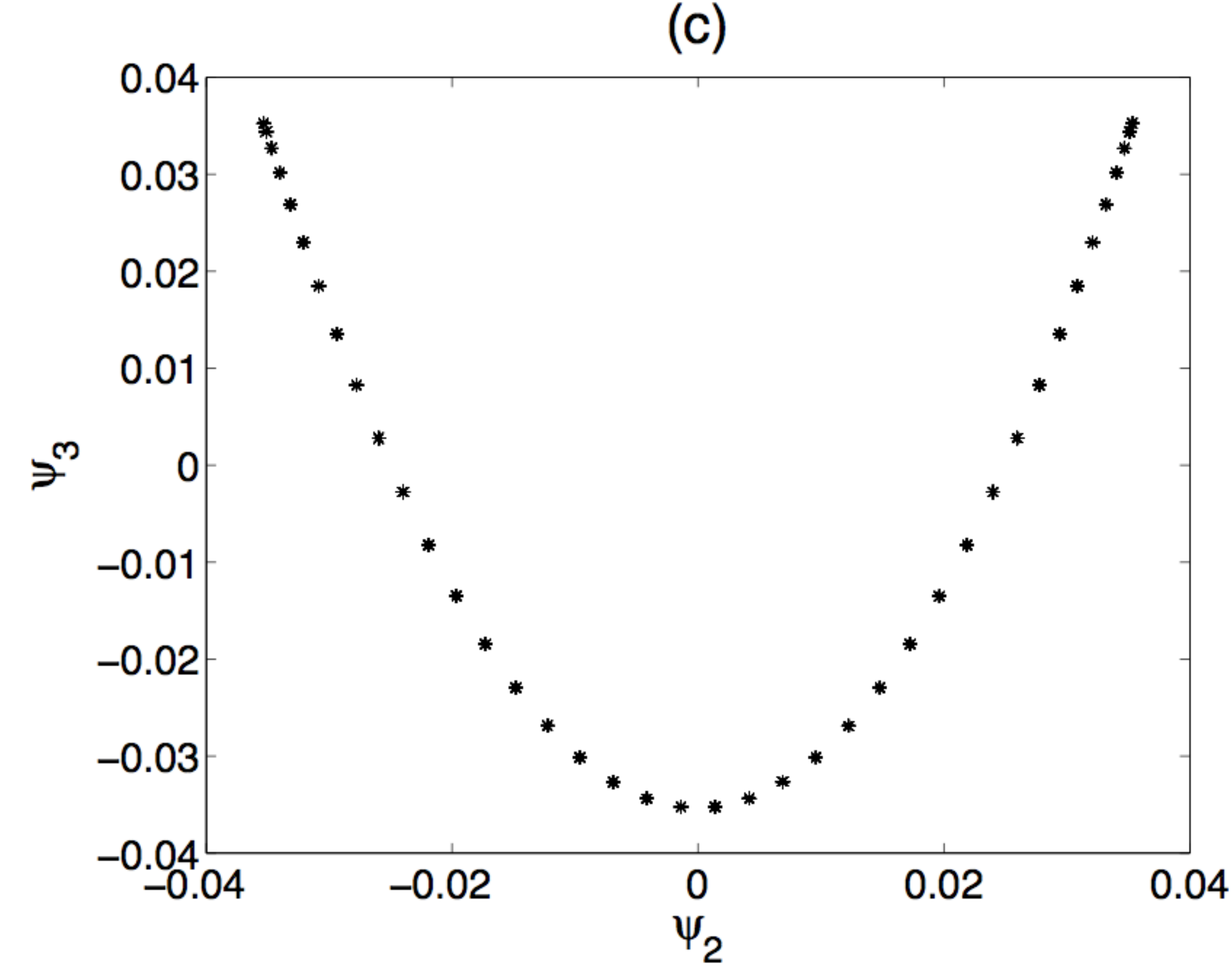}
   }
    \subfigure{
  \includegraphics[scale=0.25]{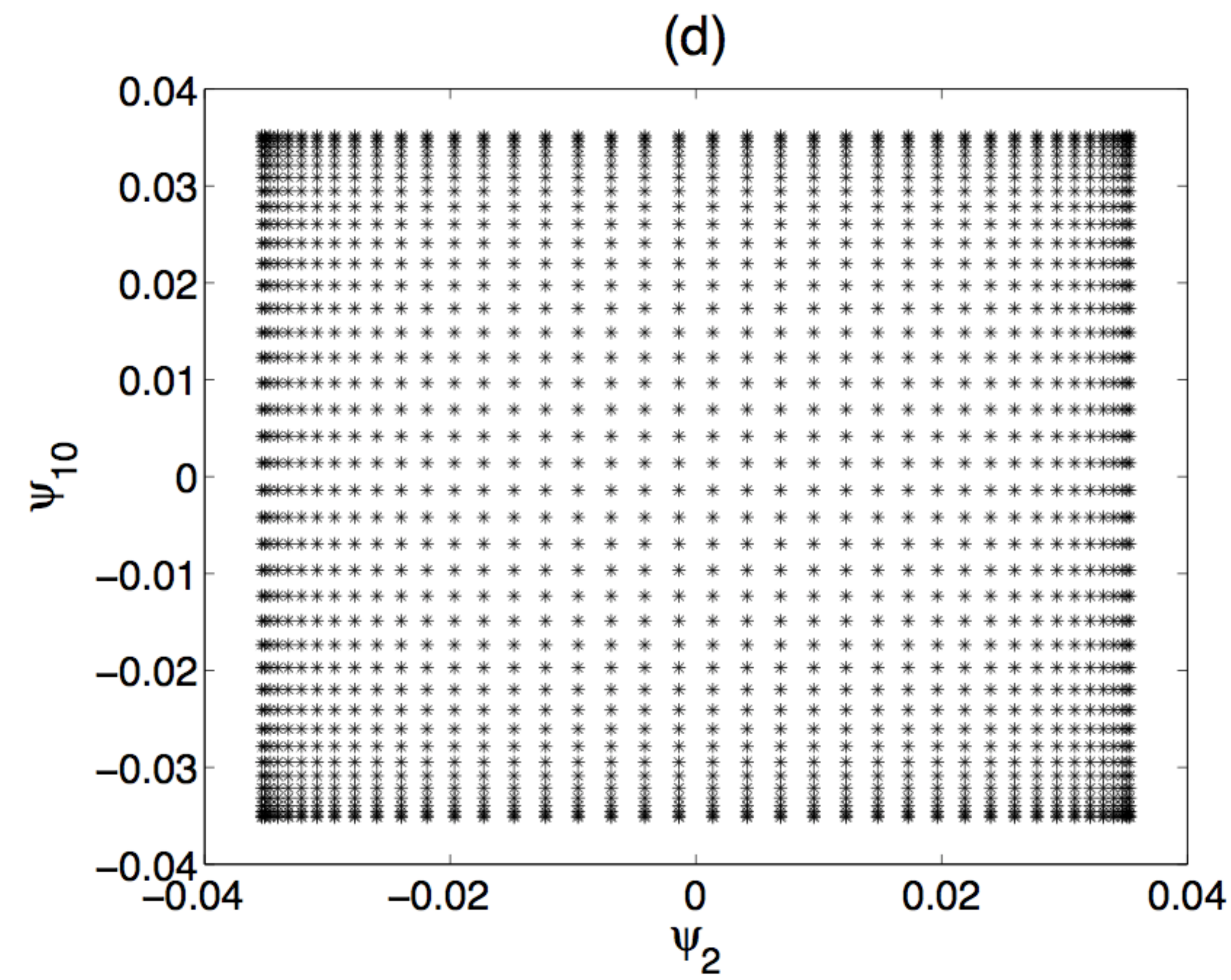}
   }
 \caption[Optional caption for list of figures]{%
The effect of data sampling. (a) A 40 by 40 array of regularly spaced points are placed on a square and subsequently wrapped around a cylinder. (b)
Entry in the first non-trivial eigenvector of the Markov matrix $K$ versus the first cylindrical coordinate $\theta$ for each data point. (c) Entry in the
second non-trivial eigenvector of $K$ versus the first one. (d) Entry in the third non-trivial eigenvector of $K$ versus the first one.}\label{CYL02}
\end{figure}

\begin{figure}[ht]
 \centering
 \subfigure{
  \includegraphics[scale=0.25]{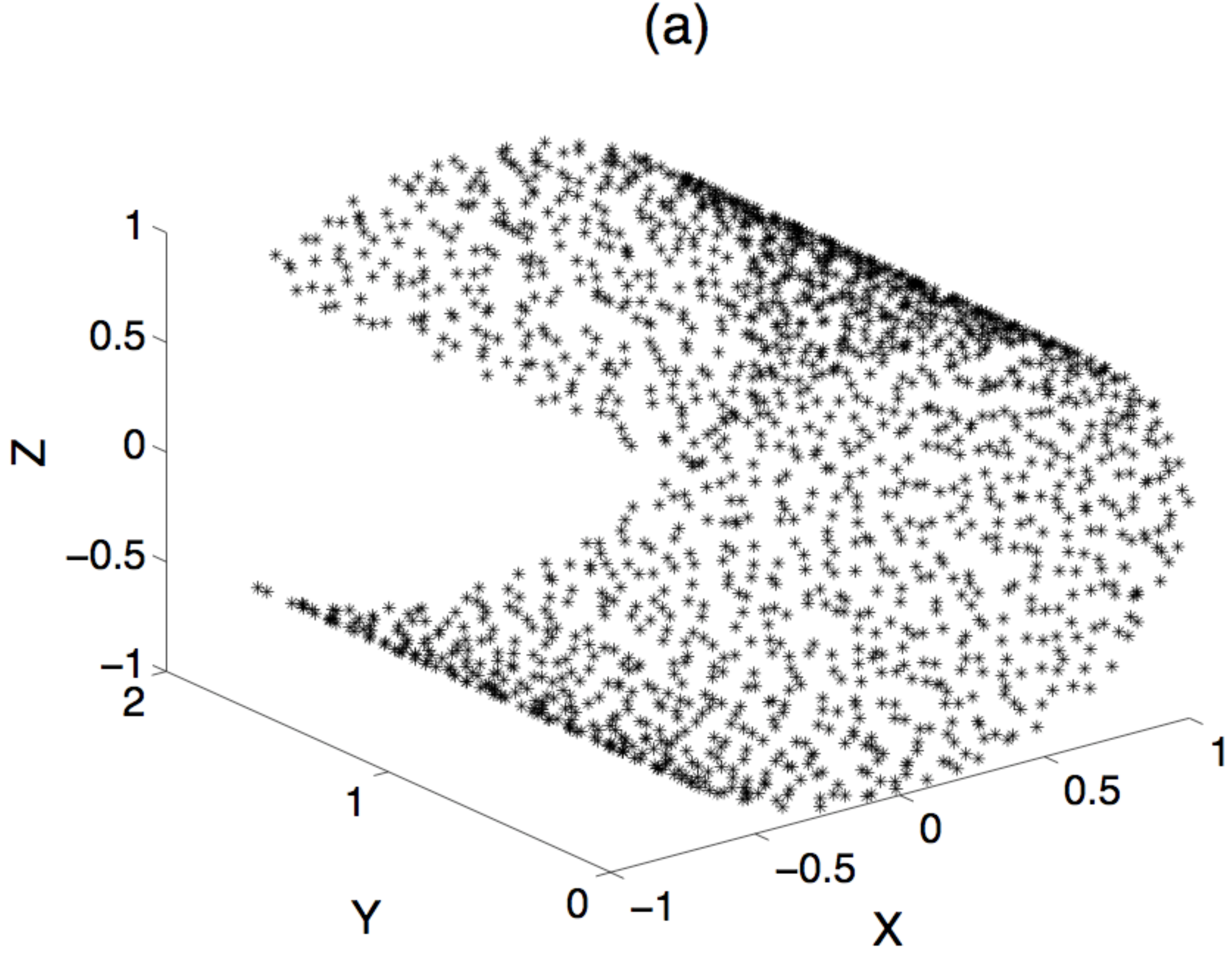}
   }
 \subfigure{
  \includegraphics[scale=0.25]{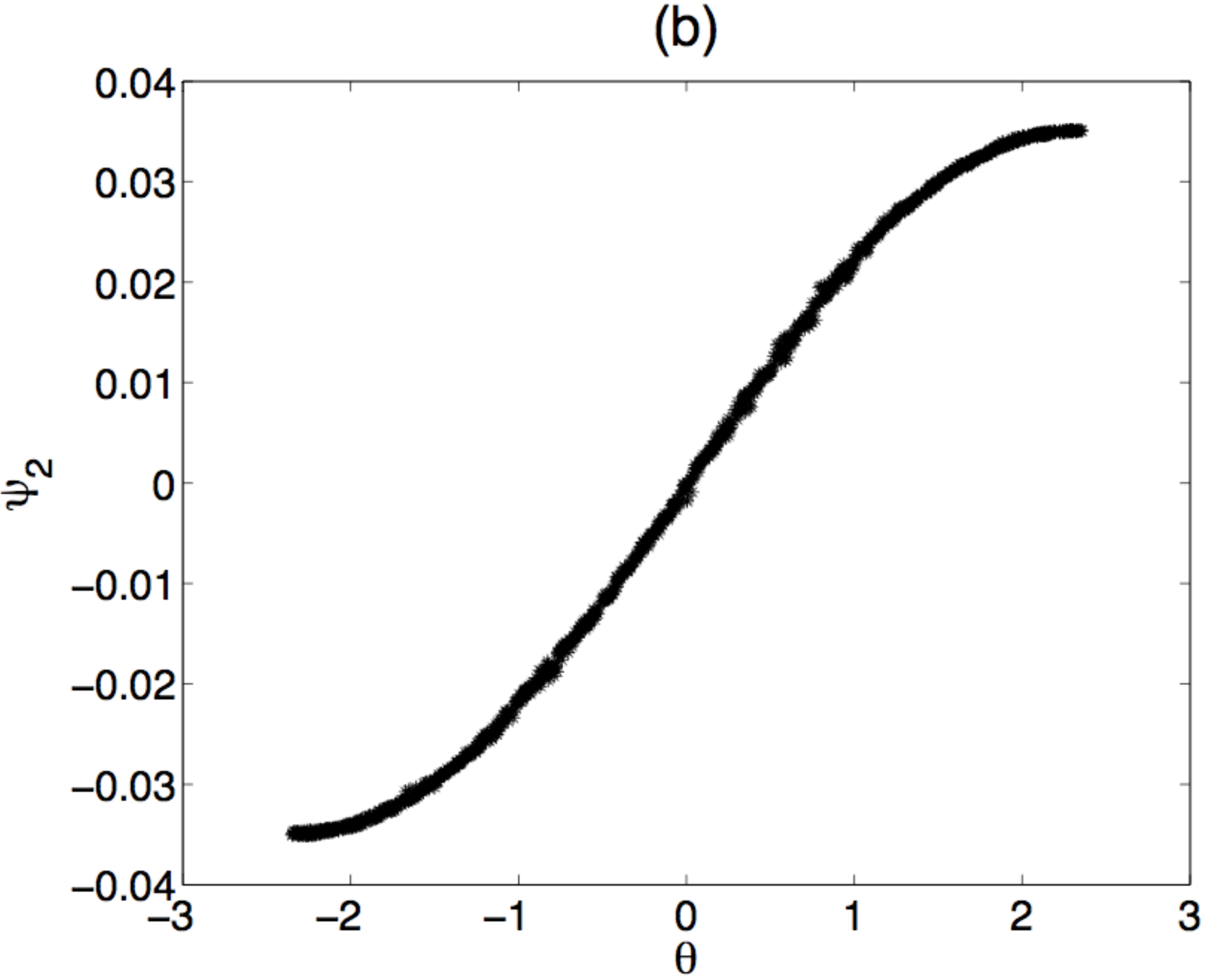}
   }
    \subfigure{
  \includegraphics[scale=0.25]{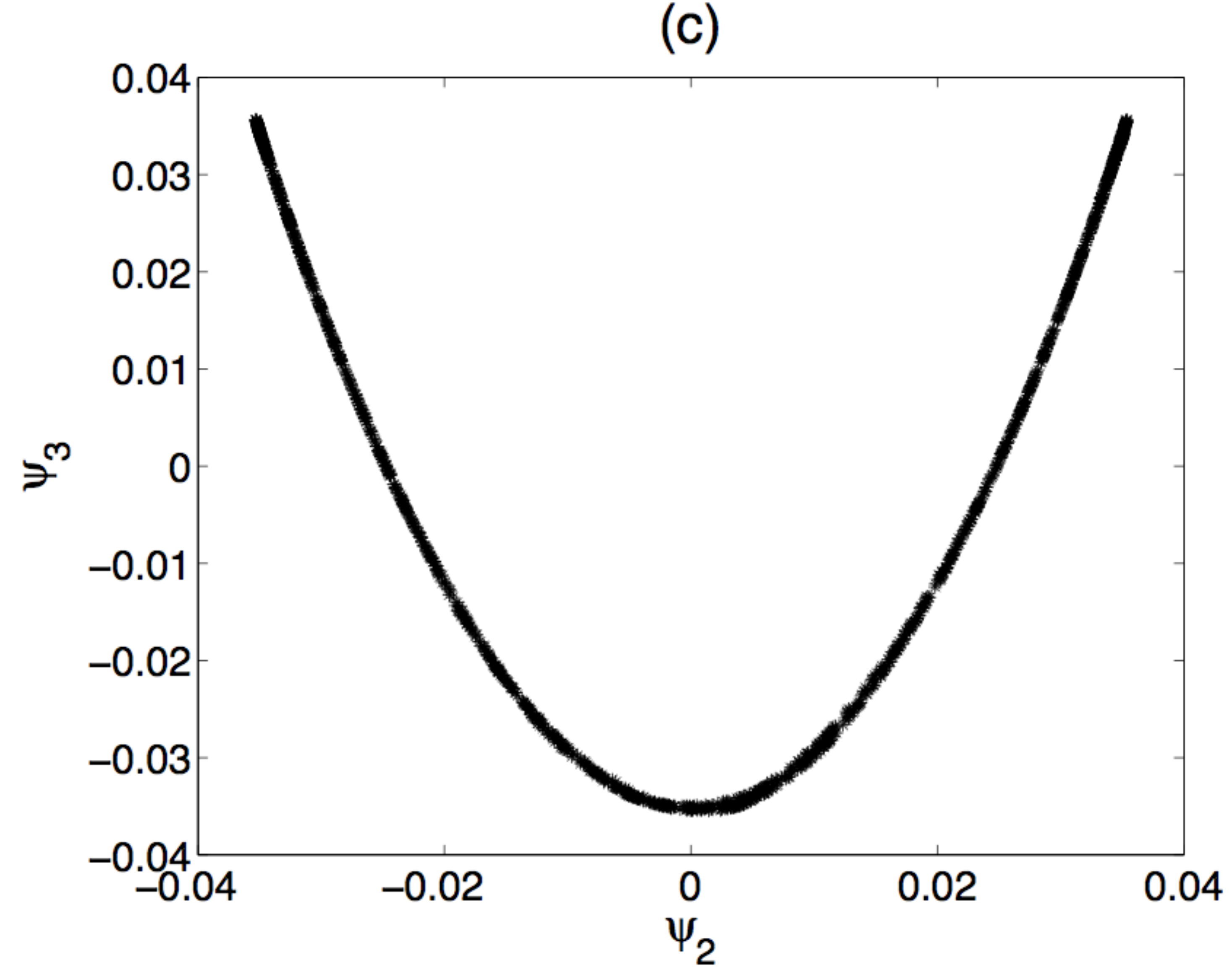}
   }
    \subfigure{
  \includegraphics[scale=0.25]{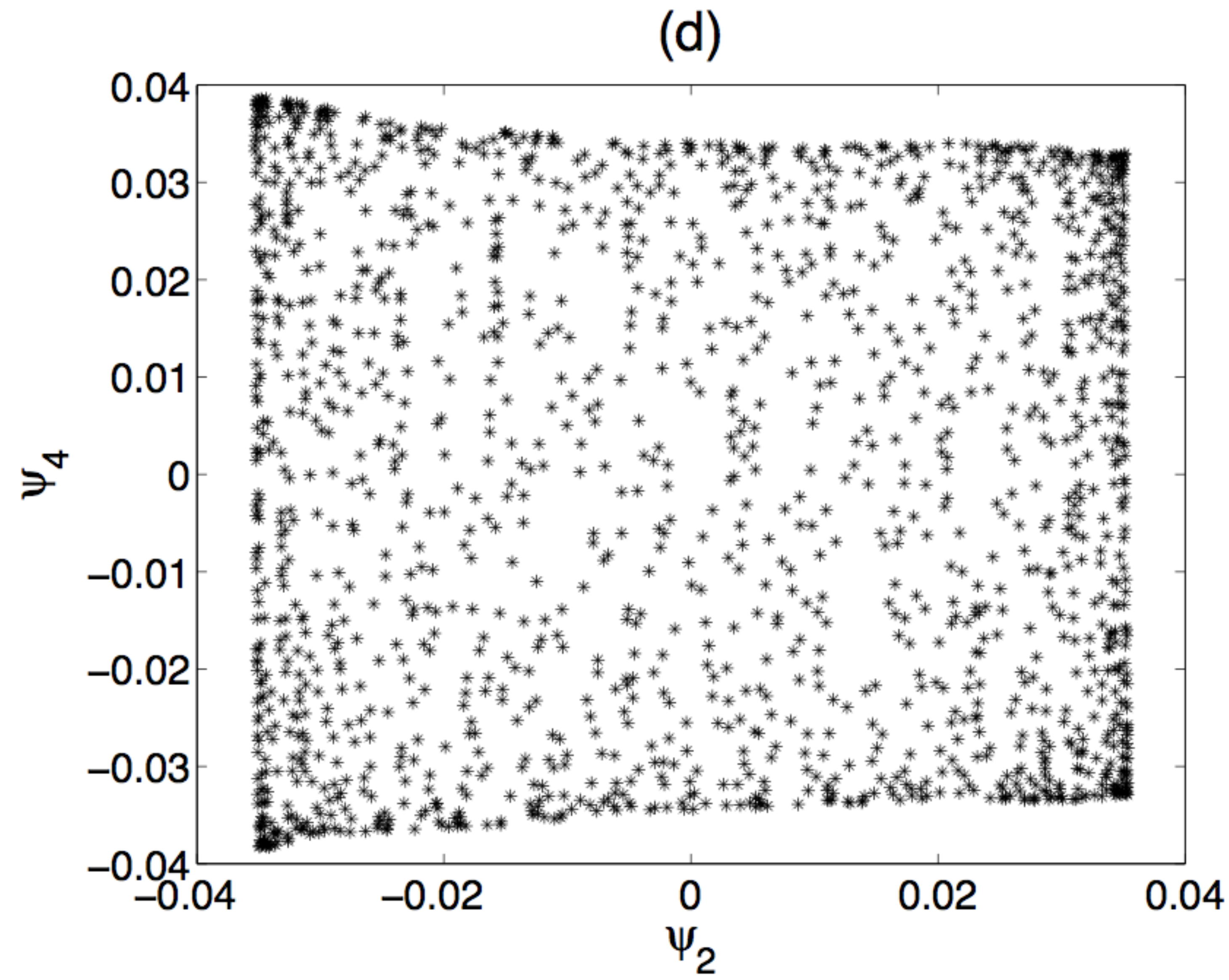}
   }
 \caption[Optional caption for list of figures]{%
More on the effect of data sampling: noise. (a) 1600 points are initially randomly placed in each of the 40 by 40 array small squares forming the unit square and afterwards bent around a cylinder. (b) Entry in the first non-trivial eigenvector of the Markov matrix $K$ versus the first cylindrical coordinate $\theta$ for each data point. (c) Entry in the second non-trivial eigenvector of $K$ versus the first one. (d) Entry in the third non-trivial eigenvector of $K$ versus the first one}\label{CYL03}
\end{figure}

%
Irregularity of sample points can be easily seen to lead to problems in this simple example.
Consider two additional cases, for different sample point distributions: First, a 40 by 40 array of regularly spaced points are placed on a square, and subsequently wrapped around the same cylinder (Fig. \ref{CYL02}(a)). Second, 1600 points are initially randomly placed in each of the 40 by 40 array small squares forming the unit square and afterwards bent around the cylinder (Fig. \ref{CYL03}(a)). As clearly visible in Figs. \ref{CYL02}(b)(c)(d) and Figs. \ref{CYL03}(b)(c)(d), this time dependencies between eigenvectors are very well defined.

\begin{figure}[t]
 \centering
 \includegraphics[scale=0.43]{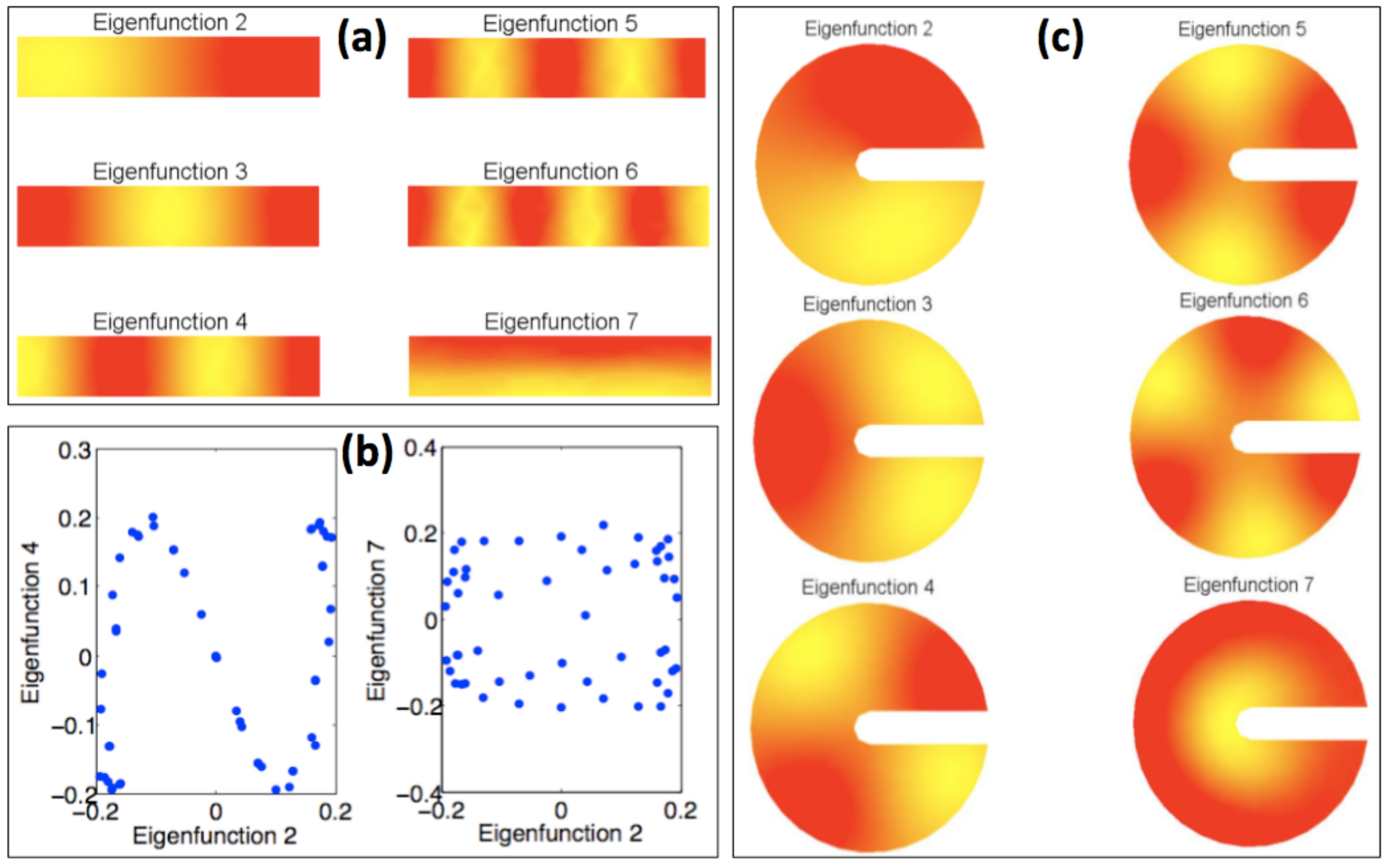}
\caption{The analogy between traditional diffusion on domains and diffusion on graphs from sampled data.
%
%
(a) The solution to the finite element method (FEM) formulation of the PDE (partial differential equation) eigenvalue problem $\nabla^2 \phi  = \lambda \phi$ with no flux boundary conditions is reported for a narrow two-dimensional rectangular stripe. The second and seventh eigenfunctions are found to be uncorrelated and suitable to parametrize the two relevant dimensions of the manifold. (b)  Entries in the first non-trivial eigenfunction of the problem in figure (a) versus entries in the fourth eigenfunction (sampled at scattered locations of the computational domain) reveals a strong correlation between those two functions: the fourth eigenvector (which we know corresponds to the third harmonic, $cos(3 \bar x)$, does not encode a new direction on the data manifold. Right-hand side: Entries in the first non-trivial eigenfunction of the problem in figure (a) versus entries in the seventh eigenfunction (at the same scattered locations) confirms that the seventh eigenvector (which we know corresponds to $cos(\bar y)$, encodes a new, second direction. (c) Different domain, same premise: The solution to the finite element method (FEM) formulation of the PDE (partial differential equation) eigenvalue problem $\nabla^2 \phi  = \lambda \phi$ with no flux boundary conditions is reported for a two-dimensional manifold with complicated boundary. The second and seventh eigenfunctions are found to be uncorrelated and suitable to parametrize the two relevant dimensions of the manifold (an ``angular" and a ``radial" one).} \label{LaplacianEig01}
\end{figure}

While the first non-trivial eigenvector $ \psi_2$ always characterizes the principal direction on the manifold, no general recipe can be formulated for an {\em a priori} identification of the subsequent uncorrelated eigenvectors parameterizing other dimensions.
We have already seen that eigenvectors in (\ref{DMAPdef}) are often dependent; this implies that they do not encode new directions along
the data manifold; in this sense, they are redundant for our embedding.
In order to obtain more insight in eigenvector dependency (and, in other words, in how diffusion is linked with manifold parametrization),
consider, as our domain of interest, a narrow two-dimensional stripe - or, in our case, data points densely sampled from it.
Fig. \ref{LaplacianEig01}(a) reports the solution to the discretized (through the finite element method, FEM) eigenvalue problem $\nabla^2 \phi  = \lambda \phi$
with Neumann boundary conditions.
The first non-trivial eigenfunction is analytically given by ${\rm cos}(\bar x)$ where $\bar x$ denotes the horizontal space direction,
and is very well approximated by the FEM numerics; the point to notice is that $cos(\bar x)$ is one-to-one with $\bar x$ between $0$ and $2\pi$;
so the first nontrivial diffusion eigenvector parameterizes one manifold direction (the $\bar x$).
Several subsequent eigenfunctions still correlate with the $\bar x$ direction: they are simply higher harmonics (${\rm cos}(2 \bar x)$, ${\rm cos}(3 \bar x)$,...).
We have to go as high as  the seventh eigenfunction (which analytically is ${\rm cos}(\bar y)$) to find something that is one-to-one with the second, independent, vertical direction $\bar y$ (see Fig.  \ref{LaplacianEig01}(b) where the first non-trivial eigenfunction is plotted against both the fourth and seventh eigenfunction at scattered locations).
A more complex two dimensional geometry is considered in Fig. \ref{LaplacianEig01}(c).
Similarly to the above example, the first non-trivial eigenfunction parameterizes one of the manifold ``principal dimensions" (the angular coordinate) , while the next (seventh) uncorrelated eigenfunction can be used to parameterize the other relevant (radial) coordinate
(it is just an accident that we had to go to {\em seventh} eigenfunction in both cases).
In practical applications, only a discrete set of sample points on the manifold in question is available as an input.
Starting from those points, the Diffusion Maps create a graph, where the points are the graph nodes and the  edges are weighted on the basis of point distances, as described above.
Noticing that the (negatively defined) normalized graph Laplacian ${L}$ is given by \cite{GraphLap}:
\begin{equation}\label{Laplacian}
L=D^{-1} W -I,
\end{equation}
with $I$ being the $M \times M$ identity matrix, we immediately recognize the link between the eigenvalue problem in Fig. \ref{LaplacianEig01} and the mapping (\ref{DMAPdef}) based on the spectrum of the Markov matrix (\ref{Markov}).
%
%
%
%
%
%
%
Diffusion on this graph (i.e. obtaining the spectrum of the graph Laplacian) approximates,
at the appropriate limit \cite{Coifman05pnas01}, the usual diffusion in the original domain;
it provides an alternative -different than our FEM, irregular mesh-  discretization of the Laplace equation eigenproblem in the original domain,
and asymptotically recovers the spectrum of the Laplace operator there.

\section{The proposed approach}\label{approach}
We demonstrate the feasibility of constructing reduced kinetics models for combustion applications, by extracting the slow dynamics on a manifold globally parameterized by a truncated diffusion map. We focus on spatially homogeneous reactive mixtures of ideal gases under fixed total enthalpy $H$ and pressure $P$. Such a set-up is relevant for building up tables to be used in reactive flow solvers in the low Mach number regime.
In such systems, a complex reaction occurs with $n$ chemical species \{$A_1$,...,$A_n$\} and $d$ elements involved in a (typically) large number, $r$ of elementary steps:
\begin{equation}\label{elementarystep}
\sum\limits_{p = 1}^n {\alpha _{sp} A_p  \mathbin{\lower.3ex\hbox{$\buildrel\textstyle\rightarrow\over
{\smash{\leftarrow}\vphantom{_{\vbox to.5ex{\vss}}}}$}} } \sum\limits_{p = 1}^n {\beta _{sp} A_p ,\quad s = 1,...,r},
\end{equation}
where $\alpha_{sp}$ and $\beta_{sp}$ represent the stoichiometric coefficients of the $p$-th species in the $s$-th step. Time evolution of chemical species can be modeled by a system of ordinary differential equations (ODEs) cast in the general form:
\begin{equation}\label{detailedODE}
\frac{{d y}}{{dt}} = \sum\limits_{s = 1}^r { \gamma _s \Omega _s \left( { y,T} \right)},
\end{equation}
with $ \gamma _s  = \left( {\beta _{s1}  - \alpha _{s1} ,...,\beta _{sn}  - \alpha _{sn} } \right)$, while the reaction rate function $\Omega_s$  is usually expressed in terms of the concentration vector $ y$ by mass action laws and Arrhenius dependence on the temperature $T$ .
Clearly, a constraint on a thermodynamic potential is required in order to close the system (\ref{detailedODE}), thus providing an additional equation for temperature. Below (without loss of generality) we consider reactions under fixed total enthalpy $H$.

The first step of our method consists in the identification of a discrete set of states lying in a neighborhood of the low-dimensional attracting manifold.
While many possible constructions have been suggested in the literature (see, e.g., \cite{ILDM,MaasGoussis,ChiavazzoKarlinPRE,Chiavazzo2012}) here, in the spirit of the {\em equation free} approach \cite{EqFree01,EqFree02}, we assume that we have no access to the analytical form of the vectorfield; instead, we only have
access to a ``black box" subroutine that evaluates the rates $ f ( y)$, and, when incorporated in a numerical initial value solver, can provide
simulation results.

\subsection{Data collection}
\begin{figure}[t]
 \centering
 \includegraphics[scale=0.8]{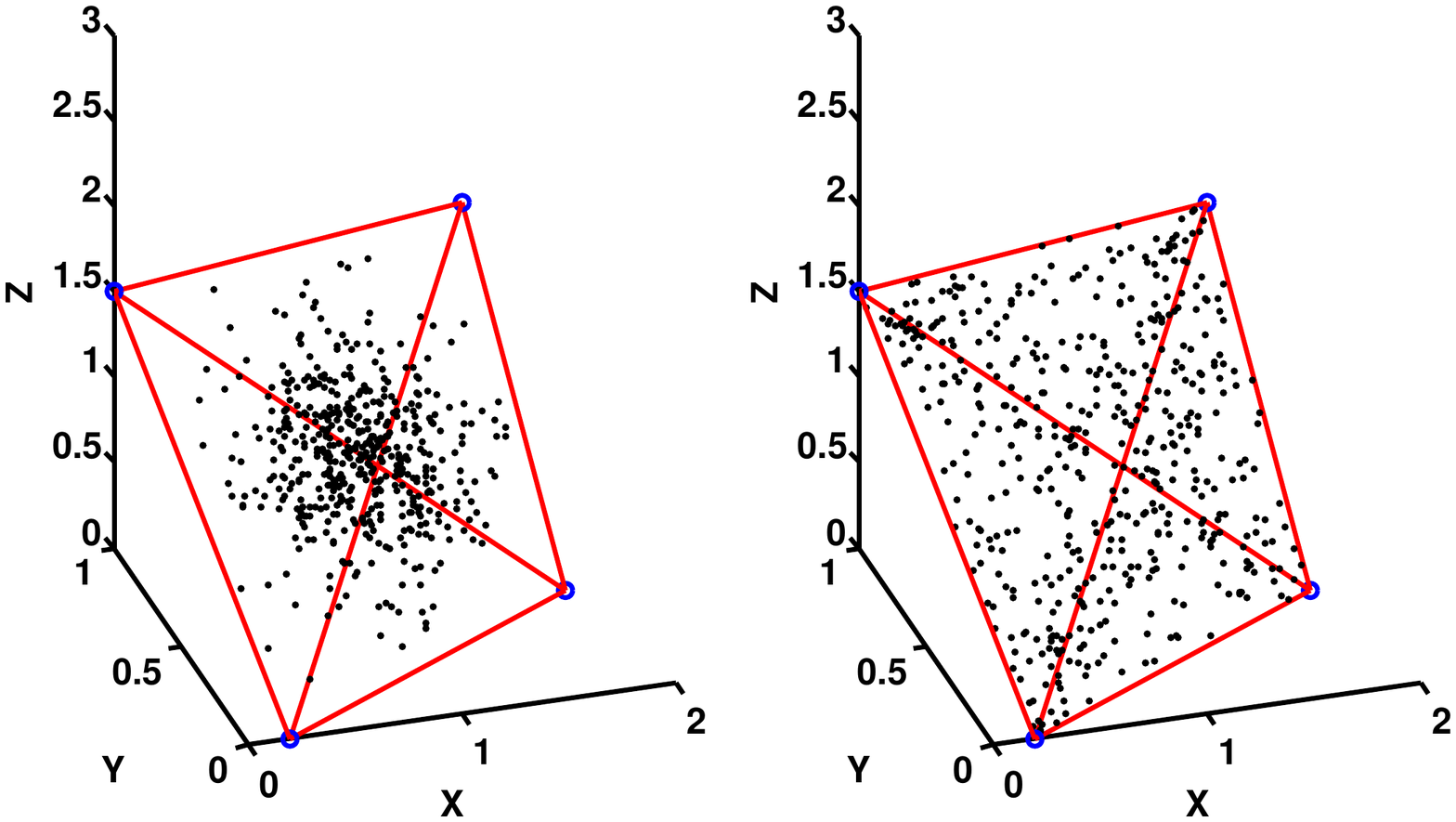}
\caption{On sampling initial conditions in a convex polytope in $\Re^3$ with vertices $A=(1.8,0.5,0)$, $B=(1,0,3)$, $C=(0,1,1.5)$ and $D=(0.2,0,0)$.
Left-hand side: Five hundred points are generated by Eq. (\ref{randomcomb}) with uniformly random values $0 \le \tilde w_i \le 1$ and $\bar w_i  = \tilde w_i /\sum\nolimits_{j = 1}^4 {\tilde w_j}$; notice the poor sampling close to the boundaries. Right-hand side: Five hundred points are generated by Eq. (\ref{randomcomb}) with uniformly random values $0 \le z_i \le 1$, $\tilde w_i  = \left[ { - \ln \left( {z_i } \right)} \right]^{1.5}$ (i.e. $p=1.5$) and $\bar w_i  = \tilde w_i /\sum\nolimits_{j = 1}^4 {\tilde w_j}$. The latter approach generates a more uniform sampling of the polytope interior.} \label{sampling}
\end{figure}

To start the procedure, we need an ensemble of representative data points on (close to) the manifold we wish to parametrize.
To ensure good sampling, our ensemble of points comes from integrating
Eqs. (\ref{detailedODE}) starting from a (rich enough) set of random states within the admissible phase-space (a convex polytope defined by elemental conservation constraints and concentration positivity).
After sufficient time to approach a neighborhood of the manifold, samples are collected from each such trajectory.
As a result, a set of points $\{{\bf y}_i, i=1,...,M\}$ in $\Re^n$ (hopefully dense enough within the region of interest) becomes available for defining the manifold.
To construct the required initial conditions we first search for all vertices of the convex polytope defined by a set of equalities and inequalities as follows:
\begin{equation}\label{polytopedef}
\begin{array}{l}
 \sum\limits_{\alpha = 1}^n { y_\alpha c_{\alpha \beta} / \bar W_\alpha  = } \sum\limits_{\alpha = 1}^n { y_\alpha^{eq} c_{\alpha \beta} / \bar W_\alpha ,} \quad \forall \beta = 1,...,d \\
  y_\alpha  > 0,\quad \forall \alpha = 1,...,n, \\
 \end{array}
\end{equation}
where $c_{\alpha \beta}$ and $\bar W_\alpha$ denote the number of atoms of the $\beta$-th element in the species $\alpha$ and the molecular weight of species $\alpha$, respectively, while the state vector $ y = \left( { y_1 ,..., y_n } \right)$ expresses species concentration in terms of mass fractions.
Selection of random initial conditions is performed by convex linear combinations of the $v$ polytope vertices $\{ y^{pol}_i \}$:
\begin{equation}\label{randomcomb}
 y^{in}=\sum_{i=1}^{v} \bar w_i  y^{pol}_i,
\end {equation}
with $\{\bar w_i\}$ being a set of $v$ random weights such as $\sum_{i=1}^{v} \bar w_i =1$. Clearly, owing to convexity, equation (\ref{randomcomb}) always provides states within the admissible space.
In combustion applications, the phase-space region of interest goes from the fresh mixture conditions to the thermodynamic equilibrium $ y^{eq}$, hence in eq. (\ref{randomcomb}) we consider a subset of the polytope vertices $\{ y^{pol}_i\}$ based on their vicinity (in the sense specified in the Appendix) to the mean point of the mixing line connecting the fresh mixture point to equilibrium.  It is worth noticing that, upon the choice of $v$ random numbers $\{\tilde w_i, i=1,...,v\}$ uniformly distributed over the range $0 \le \tilde w_i \le 1$, weights might be straightforwardly obtained by normalization: $\bar w_i  = \tilde w_i /\sum\nolimits_{j = 1}^v {\tilde w_j}$. However, such an approach leads to poor sampling in the vicinity of the polytope edges and, at the same time, to oversampling within its interior.
Therefore, in order to achieve a more uniform sampling in the whole phase-space region of interest, the weights are chosen as follows:
\begin{equation}
\tilde w_i  = \left[ { - \ln \left( {z_i } \right)} \right]^p ,\quad \bar w_i  = \tilde w_i /\sum\limits_{j = 1}^v {\tilde w_j } ,\quad i = 1,...,v
\end{equation}
with $z_i$ representing random values uniformly distributed within the interval $0 \le z_i \le 1$, and $1\le p \le 2$ a free parameter (see also Fig. \ref{sampling}).

Trajectories starting at the random initial conditions $ y^{in}$ computed by (\ref{randomcomb}) are evolved for $\tau_f$, after which samples are collected as they proceed towards the equilibrium point $ y^{eq}$.
Samples from the same trajectory are retained if their distance (see Appendix) exceeds a fixed threshold. 
%
We would like the sample to be as uniform as possible in the original space (which we will call the ambient space) because doing so yields a better parameterization with Diffusion Maps \cite{GearPreprint,SondayPrep}.
However, such a condition is not naturally fulfilled by samples of time integration: the trajectories (hence also our sampled points) often
show a tendency to gather in narrow regions (especially close to the equilibrium point, governed by the eigenvalue differences in the linearized
dynamics).
Hence, we also performed an {\em a posteriori} data filtering (subsampling) where neighbors within a minimum distance $d^{min}$ are removed.

The diffusion maps approach is performed as outlined in Section \ref{DMAP}, where distances $d_{ij}$ are computed as illustrated in the Appendix, whereas the parameter $\varepsilon$ in (\ref{similarity}) can be chosen as a multiple of the quantity $\max_j \min_{i \ne j} d_{ij}$ \cite{GearPreprint,Maggioni2011,SondayPrep}. \footnote{A better choice for $\varepsilon$ is to make it a multiple of what we will call the {\em critical diffusion distance}:  the maximum edge length such that, if all edges of at least that length are deleted in the distance graph, the graph becomes disconnected.  The reason this distance is important is that if $\varepsilon$ is much smaller than this, the diffusion map will find disjoint sets.} 
The diffusion map process provides a mapping from each point, $y_i$, in the ambient space to the reduced representation $u_i = [\psi_{i,2}, \psi_{i,3},..., \psi_{i,m+1}]^T$ in the $m$-dimensional reduced space. We will refer to this as the $u$-space.
The manifold, ${\bf \Omega}$, in the ambient ($y$) space is known only by the finite set of points $\{{\bf y}_i\}$ on ${\bf \Omega}$ and its mapping to $u$-space is known only up to the mapping of that set of points to the corresponding set $\{{\bf u}_i\}$.
Clearly we can use any interpolation technique to compute $y$ for any other value of $u$.  Let us call this $y = \Theta(u)$.
If $u$ is in a $m$-dimensional space, this mapping defines an $m$-dimensional manifold in $y$-space, ${\bf \Omega}_c$. If we chose an interpolation method such that $y_i = \Theta(u_i)$ then  ${\bf \Omega}_c$ contains the original set $\{{\bf y}_i\}$ but is an approximation to the slow manifold ${\bf \Omega}$.

We will also assume that we can construct a mapping in the other direction, $u = \psi(y)$ where $u_i = \psi(y_i)$ for all $i$. Finally, in the third step, we need to conceptually recast (\ref{detailedODE}), which has the form $dy/dt = f(y)$, into the reduced space as:
\begin{equation}
\frac{{du}}{{dt}} = g\left( u \right).
\end{equation}
In other words, given a value of $u$ we need a computational method to evaluate $g(u)$, yet all we have available is a method to compute $f(y)$.  To do this we have to execute the following three substeps:
\begin{enumerate}
\item Compute the $y$ on ${\bf \Omega}_c$ corresponding to the current $u$ (using whatever form of interpolation we chose earlier);
\item Compute $dy/dt = f(y)$;
\item Compute the equivalent $du/dt$.
\end{enumerate}
Since ${\bf \Omega}_c$ is only an approximation to ${\bf \Omega}$, it is highly unlikely that $dy/dt$ lies in the tangent plane of ${\bf \Omega}_c$ at the point $y$.  (If it did the problem of computing an equivalent $du/dt$ would be straightforward.)
Two possible solutions to this dilemma are (i) project $dy/dt$ onto the tangent plane, or (ii) extend the mapping $u = \psi(y)$ to include a neighborhood of ${\bf \Omega}_c$ (a many-to-one map). If we do the latter, we can write:
\begin{equation}\label{reducedODE}
\frac{{du}}{{dt}} = \frac{{\partial \psi }}{{\partial y}}\frac{{dy}}{{dt}}.
\end{equation}
These two approaches are really the same, since a local extension of $\psi$ to a neighborhood on ${\bf \Omega}_c$ implies a local foliation and (\ref{reducedODE}) is simply a projection along that foliation. If an orthogonal projection is used, we simply write:
\begin{equation}\label{fastfoliationP}
\frac{{du}}{{dt}} = \left( {J^T J} \right)^{ - 1} J^T \frac{{dy}}{{dt}}
\end{equation}
where $J = \frac{{\partial \Theta }}{{\partial u}}$ and $\psi(y)$ is possibly needed only for initializing (\ref{reducedODE}) in case initial conditions are available in the ambient space. 
%
\subsection{Interpolation/extension schemes}
In the following, we will review a number of possible extension (in effect, interpolation/extrapolation) schemes that might be adopted for solving the system (\ref{reducedODE}) on a learned low dimensional manifold.
\subsubsection{Nystr\"om extension.}
An established procedure for obtaining the $\alpha$-th DMAP coordinate $\psi_\alpha$ at an arbitrary state $ y \in \Re^n$ is the popular Nystr\"om extension \cite{Nyst} :
\begin{equation}\label{Nystrom}
\begin{array}{l}
 \psi_\alpha = {{\lambda _\alpha }^{-1}}\sum\limits_{i = 1}^M {k \left( {y_i,y} \right) \psi _{i,\alpha} } , \\
 k \left( {y_i , y } \right) = \left( {\sum\limits_{j = 1}^M { w \left( {y_j , y } \right)} } \right)^{ - 1} w \left( {y_i , y} \right), \\
 w \left( {y_i , y} \right) = e^{ - \left( {d_{i} /\varepsilon } \right)^2 } , \; d_i = \left\| {y'_i  - y'} \right\| ,\\
 \end{array}
\end{equation}
where $\lambda_\alpha$ and $ \left( {\psi _{1,\alpha } ,...,\psi _{M,\alpha } } \right)$ are the $\alpha$-th eigenvalue and eigenvector of the Markov matrix $K$, respectively. For the combustion case below, the $d_{i}$ denote the Euclidean distances between rescaled points as discussed in the Appendix ($  y' = R  y $, $ y'_i= R  y_i $).
%
%
%
%
The Jacobian matrix at the right-hand side of (\ref{reducedODE}) can be obtained by differentiation of (\ref{Nystrom}) as follows \cite{SondayPrep,SondayThesis}:
\begin{equation}\label{Jacobian}
\begin{array}{l}
 \frac{{\partial  \psi _\alpha  }}{{\partial  y_\beta  }} = \lambda _\alpha ^{ - 1} \left( {\sum\limits_{i = 1}^M {w \left( {y_i , y} \right)} } \right)^{ - 2} \sum\limits_{i,j = 1}^M { w \left( {y_j,y} \right)\frac{{\partial w \left( {y_i,y} \right)}}{{\partial  y_\beta  }}\left[ { \psi _{i, \alpha}   -  \psi _{j ,\alpha} } \right]}  \\
 \frac{{\partial w \left( {y_i, y } \right)}}{{\partial  y_\beta  }} = 2\varepsilon ^{ - 2} r_{\beta \beta}^2 w \left( {y_i , y } \right) \left( { y_{i,\beta}   -  y_{\beta} } \right), \\
 \end{array}
\end{equation}
where, in case of point rescaling, $r_{\beta \beta}$ is computed as indicated in the Appendix, otherwise $r_{\beta \beta}=1, \forall \beta$.
The Nystr\"om estension can be utilized for implementing the restriction operator, as well as for computing its Jacobian matrix.

\subsubsection{Radial basis functions.}
Both lifting and restriction operators may be also obtained by local interpolation through radial basis functions.
Let $u$ be a new state in the reduced space; the corresponding point in the full space $ y = \Theta \left( u \right)$ can be generally expressed by the following summation:
\begin{equation}\label{inverseRBFgeneral}
 y_\beta  = \sum\limits_{i = 1}^{nn} {\alpha_{i,\beta}  \bar \phi \left( {\left\| { u -  u_i} \right\|} \right),\quad \beta  = 1,...,n},
\end{equation}
over the $nn$ nearest neighbors of $u$ with the radial function ${\bar \phi}(\bullet)$ only depending on a distance $\left\| {\bullet} \right\|$. In this work, we focus on the following special form of (\ref{inverseRBFgeneral}):
\begin{equation}\label{inverseRBF}
 y_\beta  = \sum\limits_{i = 1}^{nn} {\alpha_{i,\beta} \left\| { u  -  u_i } \right\|^p }, \quad \beta=1,...,n
\end{equation}
where $p$ is an odd integer while $\left\| {\bullet} \right\|$ denotes the usual Euclidean distance in the reduced space.
The coefficients $\alpha_{i,\beta}$ are computed as:
\begin{equation}
\left[ {\alpha_{1,\beta}  ,...,\alpha _{nn,\beta} } \right]^T  = \Lambda ^{ - 1} \left[ { y_{1,\beta}  ,..., y_{nn,\beta}  } \right]^T ,\quad \Lambda \left( {i,j} \right) = \left\| { u_i  -  u_j } \right\|^p ,\quad \;i,j = 1,...,nn.
\end{equation}
Similarly, restriction can be expressed in the form:
\begin{equation}\label{directRBF}
 \psi _\beta  = \sum\limits_{i = 1}^{nn} {\alpha_{i,\beta}  \left\| { y' -  y'_i } \right\|^p ,\quad j = 1,...,k} ,
\end{equation}
where data in the full space have been possibly rescaled as discussed in the Appendix ($  y' = R  y $, $ y'_i= R  y_i $). The Jacobian matrix at the right-hand side of (\ref{reducedODE}) can be obtained by differentiation of (\ref{directRBF}) as follows:
\begin{equation}\label{JacobianRBF}
\frac{{\partial  \psi _\beta }}{{\partial  y_\gamma }} = p r_{\gamma \gamma}^2   \sum\limits_{i = 1}^{nn} {\alpha_{i,\beta}  \left[ {\sum\limits_{\omega = 1}^n {\left( { y'_\omega  -  y'_{i,\omega} } \right)^2 } } \right]^{\frac{p}{2} - 1} \left[ {y_\gamma -  y_{i,\gamma} }   \right]}.
\end{equation}
\subsubsection{Kriging.}
Kriging typically refers to a number of sophisticated interpolation techniques originally developed for geostatistics applications. Provided a function $f$ known on scattered data, its extension to a new point is performed via a weighted linear combination of the values of $f$ at known locations.
A noticeable feature of Kriging is that weights may depend on both distances {\em and correlations} between the available samples. In fact, one possible disadvantage of schemes only based on the quantities $\left\|  \bullet  \right\|$ (e.g. radial basis functions) is that samples at a given distance from the location where an estimate is needed are all equally treated. In contrast, Kriging offers the possibility of performing a weighting which accounts for redundancy (i.e. sample clustering) and even sample orientation.
This is done by choosing an analytical model that best fits the experimental semivariogram of the data set.
More details on Kriging can be found in Ref. \cite{KrigBook}. In this work, both interpolated points and derivatives are computed by the readily available Matlab toolbox DACE \cite{KrigSoft}.
\subsubsection{Laplacian Pyramids.}
Laplacian Pyramids are a multi-scale extension algorithm, where a function only known at $M$ (scattered) sample points can be estimated at a new location. Based on a chosen kernel and pair-wise distances between samples, this algorithm aims at generating a sequence of approximations with different resolutions at each subsequent level $l$ \cite{LP}. Let $ y $ be a new point in the full space, the $\alpha$-th coordinate of the corresponding state in the reduced space $u$ is evaluated in a multi-scale fashion as follows: $ \psi_\alpha \approx s_\alpha^{(0)}+s_\alpha^{(1)}+s_\alpha^{(2)}+...$, with
\begin{equation}\label{scalesLP}
\begin{array}{l}
 s_\alpha^{(0)}  = \sum\limits_{i = 1}^M {k^{(0)} \left( {y_i, y} \right)}  \psi_{i,\alpha} \quad \quad {\rm{for}}\;{\rm{level}}\;\;l = 0 \\
 s_\alpha^{(l)}  = \sum\limits_{i = 1}^M {k^{(l)} \left( {y_i, y} \right)}  d_{i,\alpha}^{(l)} \quad \quad \; {\rm{otherwise}}{\rm{,}} \\
 \end{array}
\end{equation}
and the differences
\begin{equation}\label{differencesLP}
\begin{array}{l}
  d_\alpha^{(1)}  =  \psi_\alpha  - s_\alpha^{(0)} \quad \quad {\rm{for}}\; \;{\rm{level}}\;\;l = 1 \\
  d_\alpha^{(l)}   =  \psi_\alpha - \sum\limits_{i = 0}^{l - 1} {s_\alpha^{(i)} } \quad {\rm{otherwise}} \\
 \end{array}
\end{equation}
are updated at each level $l$. The functions $k^l$ in (\ref{scalesLP}) are:
\begin{equation}\label{kernelLP}
\begin{array}{l}
 k^{(l)} \left( {y_i,y} \right) = q_l^{ - 1} w^{(l)} \left( {y_i,y} \right) \\
 q_l = \sum\limits_{j } {w^{(l)} \left( {y_j, y} \right)}  \\
 w^{(l)} \left( {y_i,y} \right) = \exp \left[ { - \left\| { y'_i -  y'} \right\|^2 /\sigma _l } \right] .\\
 \end{array}
\end{equation}
In Eqs. (\ref{kernelLP}), a Gaussian kernel is chosen where the parameter $\sigma_l = \sigma_0/2^l$ decreases with the level $l$, $\sigma_0$ is the fixed coarsest scale, while $y'_i$ and $y'$ denote the rescaled states as specified in the Appendix ($  y' = R  y $, $ y'_i= R  y_i $). A maximum admissible error can be set {\em a priori}, and the values $s_\alpha^{(l)}$ are only computed up to the finest level where: $\left\| { \psi _\alpha  - \sum\nolimits_k^{} {s_\alpha^{(k)} } } \right\| < {\rm{err}}$.
The use of Laplacian Pyramids for constructing a lifting operator, $ y_\alpha \approx s_\alpha^{(0)}  + s_\alpha^{(1)}  + s_\alpha^{(2)}  + ...$, is straightforward and only requires the substitution of $ \psi_{i,\alpha}$ with $ y_{i,\alpha}$ in (\ref{scalesLP}) and (\ref{differencesLP}), while Euclidean distances in the reduced space are adopted for the kernel in (\ref{kernelLP}). Based on the resemblance of Eqs. (\ref{scalesLP}) with the Nystr\"om extension (\ref{Nystrom}), it follows that:
\begin{equation}\label{JacobianLP}
\begin{array}{l}
 \frac{{\partial s_\alpha ^{(0)} }}{{\partial  y_\beta  }} = \left( {\sum\limits_{i = 1}^M {w^{(0)} \left( {y_i,y} \right)} } \right)^{ - 2} \sum\limits_{i,j = 1}^M {w^{(0)} \left( {y_j,y} \right)\frac{{\partial w^{(0)} \left( {y_i,y} \right)}}{{\partial  y_\beta  }}\left[ { \psi _{i,\alpha}  -  \psi _{j,\alpha}  } \right]}  \\
 \frac{{\partial s_\alpha ^{(1)} }}{{\partial  y_\beta  }} = \left( {\sum\limits_{i = 1}^M {w^{(1)} \left( {y_i,y} \right)} } \right)^{ - 2} \sum\limits_{i,j = 1}^M {w^{(1)} \left( {y_j,y} \right)\frac{{\partial w^{(1)} \left( {y_i,y} \right)}}{{\partial  y_\beta  }}\left[ { \psi_{i,\alpha}  -  \psi_{j,\alpha} } \right]}  - \frac{{\partial s_\alpha ^{(0)} }}{{\partial  y_\beta  }} \\
 \frac{{\partial s_\alpha ^{(2)} }}{{\partial  y_\beta  }} = \left( {\sum\limits_{i = 1}^M {w^{(2)} \left( {y_i,y} \right)} } \right)^{ - 2} \sum\limits_{i,j = 1}^M {w^{(2)} \left( {y_j,y} \right)\frac{{\partial w^{(2)} \left( {y_i,y} \right)}}{{\partial  y_\beta  }}\left[ { \psi_{i,\alpha}  -  \psi_{j,\alpha}  } \right]}  - \frac{{\partial s_\alpha ^0 }}{{\partial  y_\beta  }} - \frac{{\partial s_\alpha ^1 }}{{\partial  y_\beta  }} \\
  \vdots  \\
 \end{array}
\end{equation}
with
\begin{equation}
\frac{{\partial w^{(l)} \left( {y_i,y} \right)}}{{\partial  y_\beta  }} = 2\sigma_l^{ - 2} r_{\beta \beta}^{2} w^{(l)} \left( {y_i, y } \right)\left( { y_{i,\beta}   -  y_\beta  } \right),
\end{equation}
and the Jacobian at the right-hand side of (\ref{reducedODE}) given by
\begin{equation}\label{JacobianLP01}
\frac{{\partial  \psi _\alpha  }}{{\partial  y_\beta  }} = \sum\limits_l {\frac{{\partial s_\alpha ^{(l)} }}{{\partial  y_\beta  }}} .
\end{equation}

Similarly to RBF, LP can be applied to a subset of the sample points where, in the above procedure, only $nn$ nearest neighbors of the state $y$ ($u$) are considered for restriction (lifting).

\begin{figure}[t]
 \centering
 \subfigure{
  \includegraphics[scale=0.26]{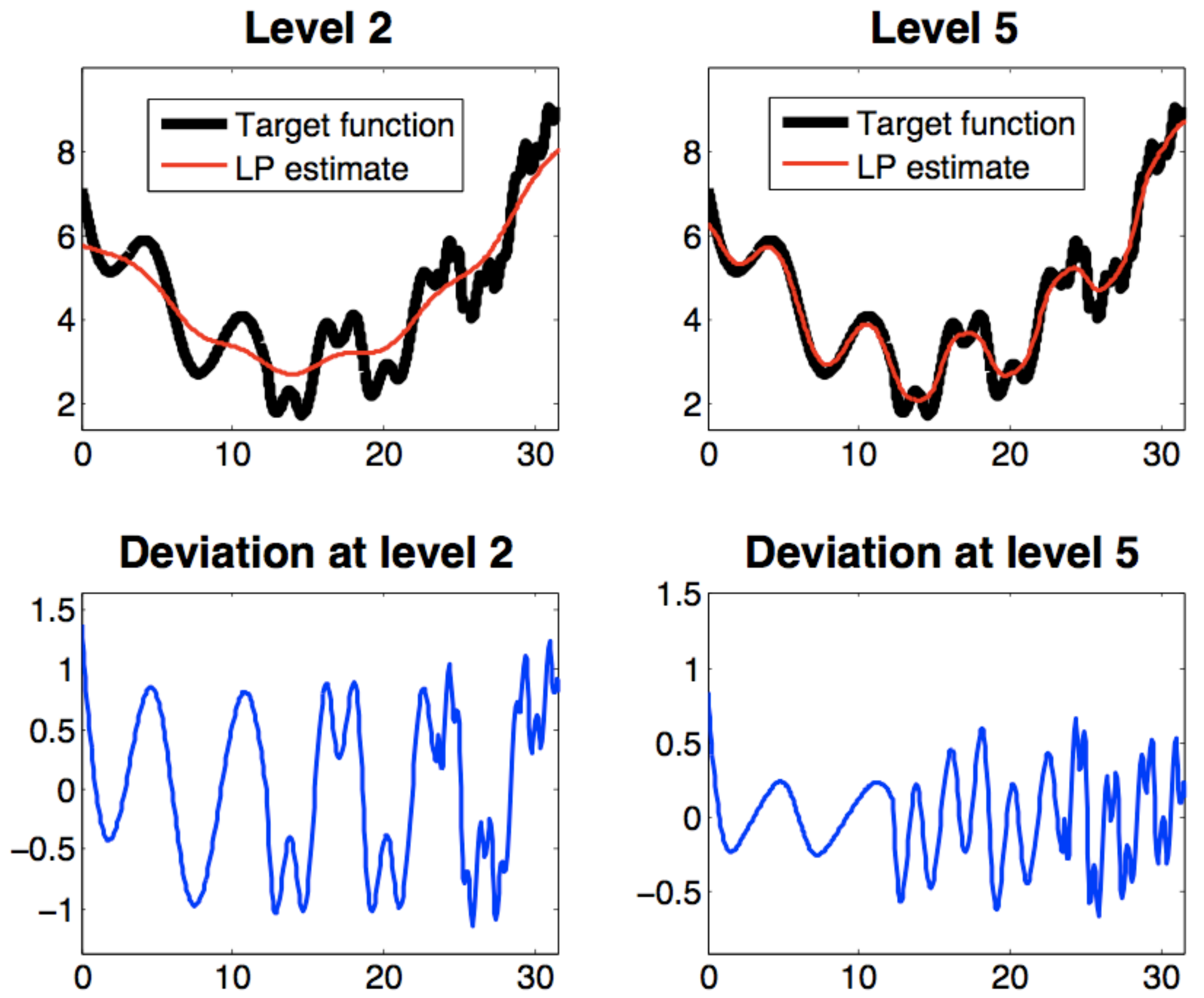}
  }
 \subfigure{
  \includegraphics[scale=0.255]{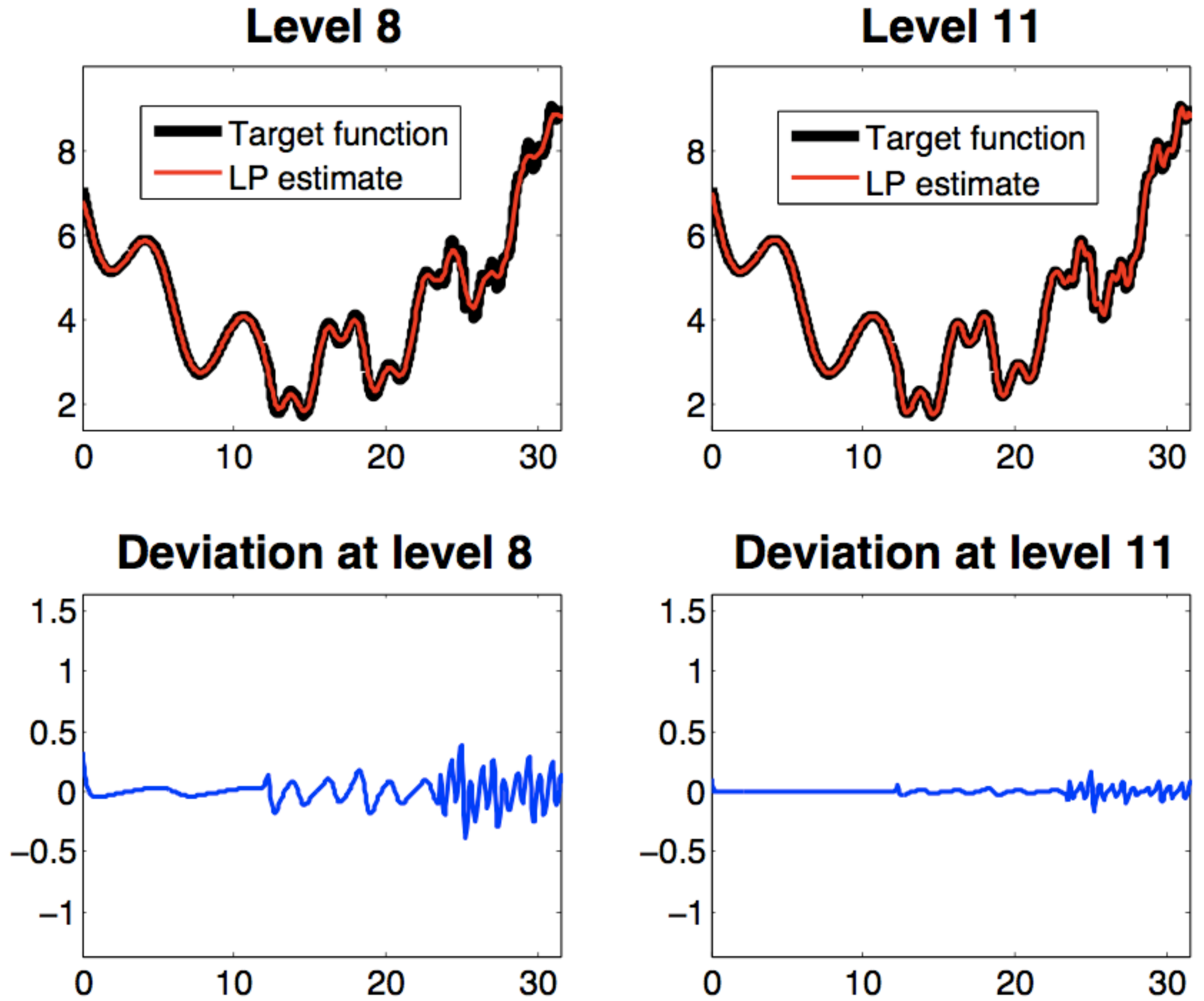}
   }
  \caption[Optional caption for list of figures]{%
Illustrating Laplacian Pyramids for a Multiscale Target Function (see text). The sample data-set is formed by 2000 points evenly distributed in the interval $[0,10 \pi]$. Top: Laplacian Pyramids used as an interpolation procedure at levels 2, 5, 8 and 11 with $\sigma_0=30$. Bottom: Difference between the true function values and the Laplacian Pyramids estimates.}\label{LP01new}
\end{figure}

\begin{figure}[t]
 \centering
 \includegraphics[scale=0.65]{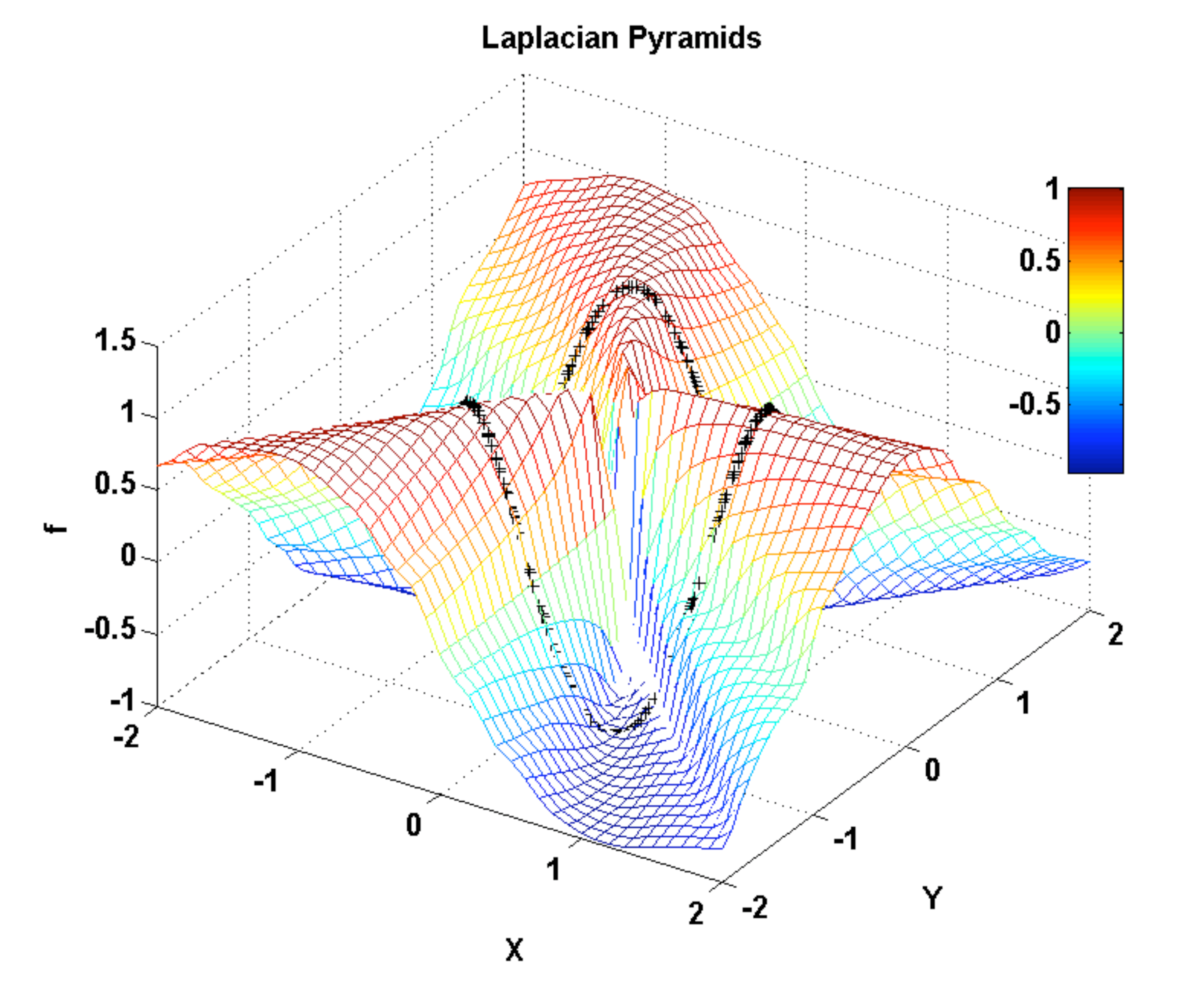}
\caption{The function $f=\cos(3 \vartheta)$, with $\vartheta=\arctan(Y/X)$, is extended to the plane $(X,Y)$ by Laplacian Pyramids (with the finest level $l=10$ and $\sigma_0=10$). The sample set $\{{\bf y}_1,...,{\bf y}_M\}$ is given by $M=350$ points randomly selected on the unit circle (black symbols).} \label{LP10}
\end{figure}

A brief explanatory illustration of the use of Laplacian Pyramids for interpolating a multi-scale function at four different levels of accuracy is
given in Figs. \ref{LP01new};  in Fig. \ref{LP10} the same scheme provides an extension of the function $f(\vartheta)=\cos(3 \vartheta)$, defined on the circle in $\Re^2$ given by $X^2+Y^2=1$ with $\vartheta=\arctan(Y/X)$.
%
%
%
%
%
%
%
%
%
\subsubsection{Geometric harmonics.}
This is an alternative multi-scale scheme for extending functions only available at $M$ scattered locations, inspired by the Nystr\"om method,
and making use of a kernel $w$ \cite{GH}.
Let $W$ be the symmetric $(M \times M)$ matrix, whose generic element reads as:
\begin{equation}\label{kernelGH}
w\left( {i,j} \right) = \exp \left[ { - d_{ij}^2 /\varepsilon _0 } \right],\quad i,j = 1,...,M
\end{equation}
with $\{ \phi_{\alpha=1,...,M}\}$ being its full set of orthonormal eigenvectors sorted according to descending eigenvalues $\{ \lambda_{\alpha=1,...,M} \}$. For $\delta>0$, let us define the set of indices $S_\delta=\{\alpha \; \;  {\rm such \; that} \; \lambda_\alpha \ge \delta \lambda_0 \}$. The extension of a function $f$ defined only at some sample points in $Z \subset \bar Z$ to an arbitrary new point in $\bar Z$ is accomplished by the following projection step (depending on the purpose, $\bar Z$ can be either the ambient space $y$ or the reduced one $u$):
\begin{equation}
f \to {\rm{P}}_\delta  f = \sum\limits_{\alpha \in S_\delta  } {\left\langle {f, \phi _\alpha } \right\rangle  \phi _\alpha } ,
\end{equation}
and the subsequent extension of ${\rm{P}}_\delta  f $:
\begin{equation}
{\bf E} f  = \sum\limits_{\alpha \in S_\delta  } {\left\langle {f, \phi _\alpha } \right\rangle \Psi _\alpha } ,
\end{equation}
where $\left\langle {\bullet,\bullet} \right\rangle $ denotes the inner product, while $\Psi _j $ reads:
\begin{equation}
\Psi _\alpha  = \lambda _\alpha^{ - 1} \sum\limits_{i = 1}^M {w\left( {y_i , y} \right) \phi _{i,\alpha} } ,\quad \alpha \in S_\delta.
\end{equation}
The above is only the first step of a multi-scale scheme, where the function $f$ is initially projected at a coarse scale with a large value of the parameter $\varepsilon_0$ in (\ref{kernelGH}). Afterwards, the residual $f- {\rm P}_\delta f$ in the initial coarse projection is projected at a finer scale $\varepsilon_1$, and so forth at even finer scales $\varepsilon_{l>1}$. A typical approach is to fix $\varepsilon_0$, and then project with $\varepsilon_l=2^{1-l}\varepsilon_0$ at each subsequent step $l$ till a norm of the residual $f- {\rm P}_\delta f$ remains larger than a fixed admissible error. Clearly, both our restriction and lifting operators can be based on Geometric Harmonics.

Similarly to RBF and LP, GH can be applied to a subset of the sample points where, in the above procedure, only $nn$ nearest neighbors of the state $y$ (or $u$) are considered for restriction (or lifting).

\begin{figure}[t]
 \centering
 \includegraphics[scale=0.5]{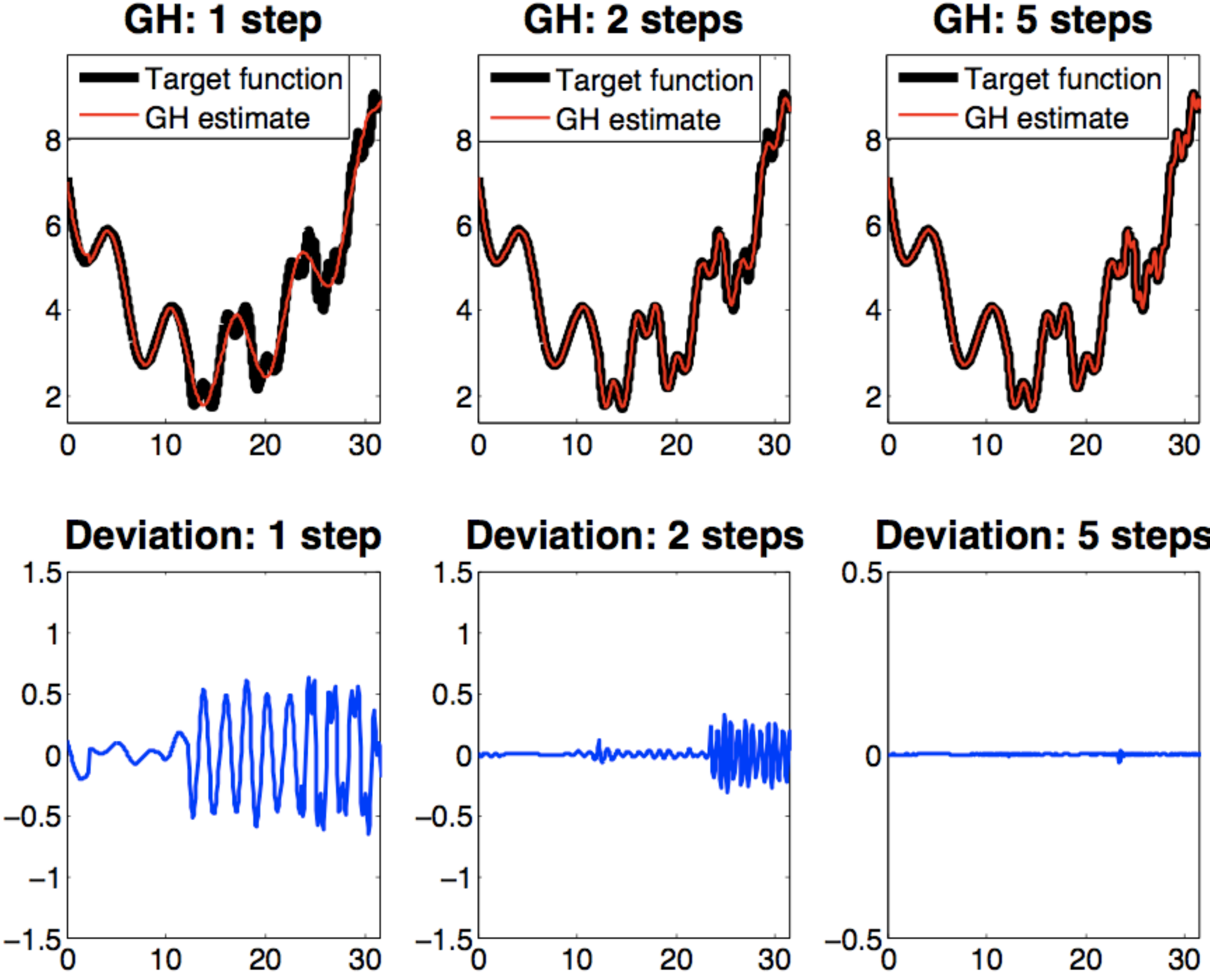}
\caption{Geometric Harmonics on a Multiscale Target Function (see text). The sample data-set is formed by 2000 points evenly distributed in the interval $[0,10 \pi]$. Top: The Geometric Harmonics (GH) scheme is used as an interpolation procedure with $\varepsilon_0=3$. Bottom: Difference between the true function values and GH estimates. From left to right: Results corresponding to 1, 2 and 8 steps are reported.} \label{GHcomparison}
\end{figure}
\begin{figure}[t]
 \centering
 \includegraphics[scale=0.85]{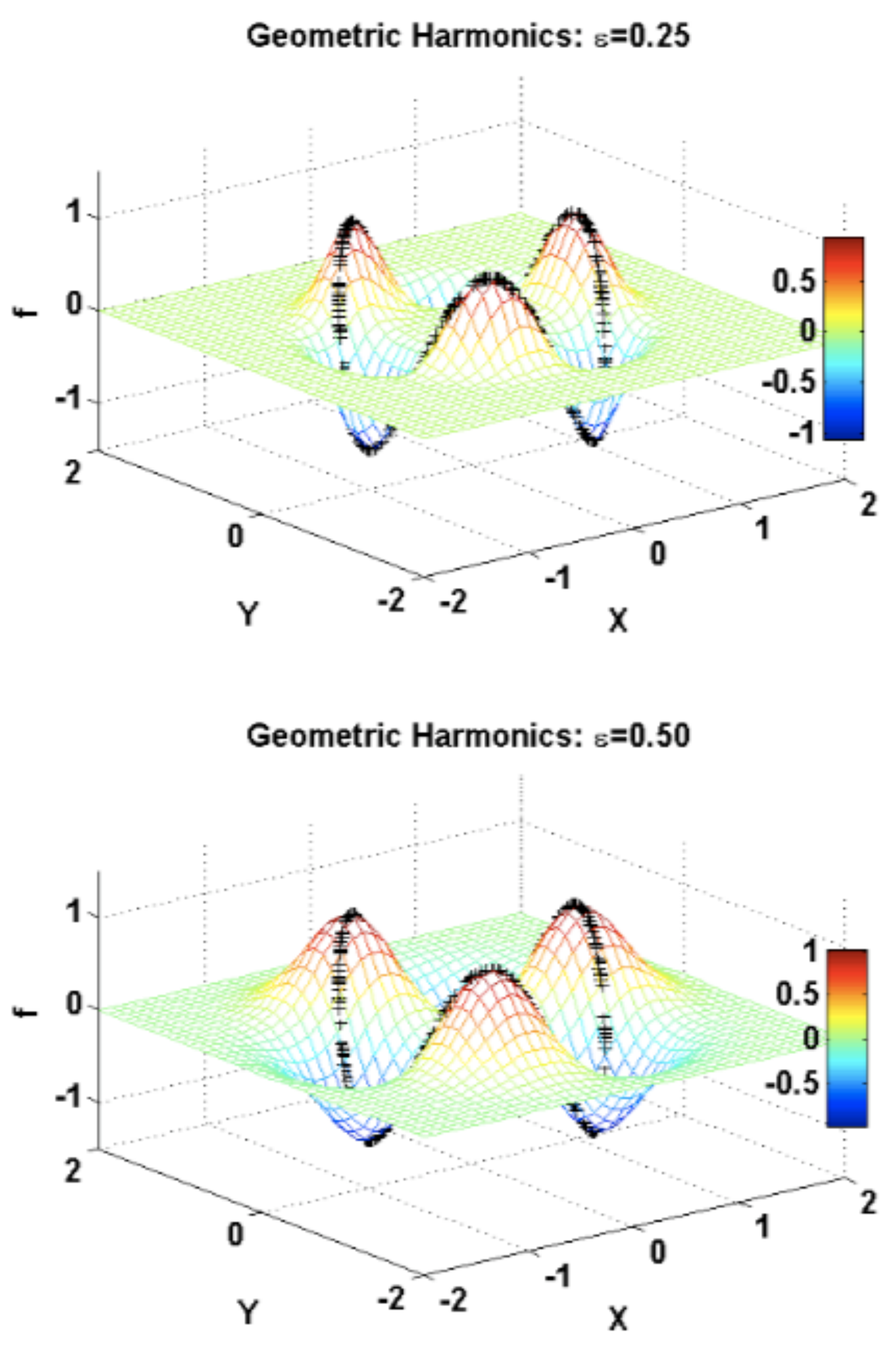}
\caption{The function $f=\cos(3 \vartheta)$, with $\vartheta=\arctan(Y/X)$, is extended to the plane $(X,Y)$ by Geometric Harmonics using $\varepsilon_0=0.25$ (top) and $\varepsilon_0=0.5$ (bottom). The sample set $\{{\bf y}_1,...,{\bf y}_M\}$ is given by $M=350$ points randomly selected on the unit circle (black symbols).} \label{GH01}
\end{figure}
Figure \ref{GHcomparison} provides an illustrative multi-scale example where the Geometric Harmonics approach is used for interpolation purposes for the
same multiscale function used in Fig. \ref{LP01new}. As expected, in the region with low-frequency components, a few steps are sufficient for accurately describing the true function, whereas more iterations are required in the  high frequency domain.
We also illustrate the use of Geometric Harmonics in extending the function $f(\vartheta)=\cos(3 \vartheta)$,  defined on the circle in $\Re^2$ given by $X^2+Y^2=1$ with $\vartheta=\arctan(Y/X)$ in Fig. \ref{GH01}.

\section{Application to an illustrative example: Homogeneous combustion}\label{results}
\begin{figure}[t]
 \centering
 \includegraphics[scale=0.7]{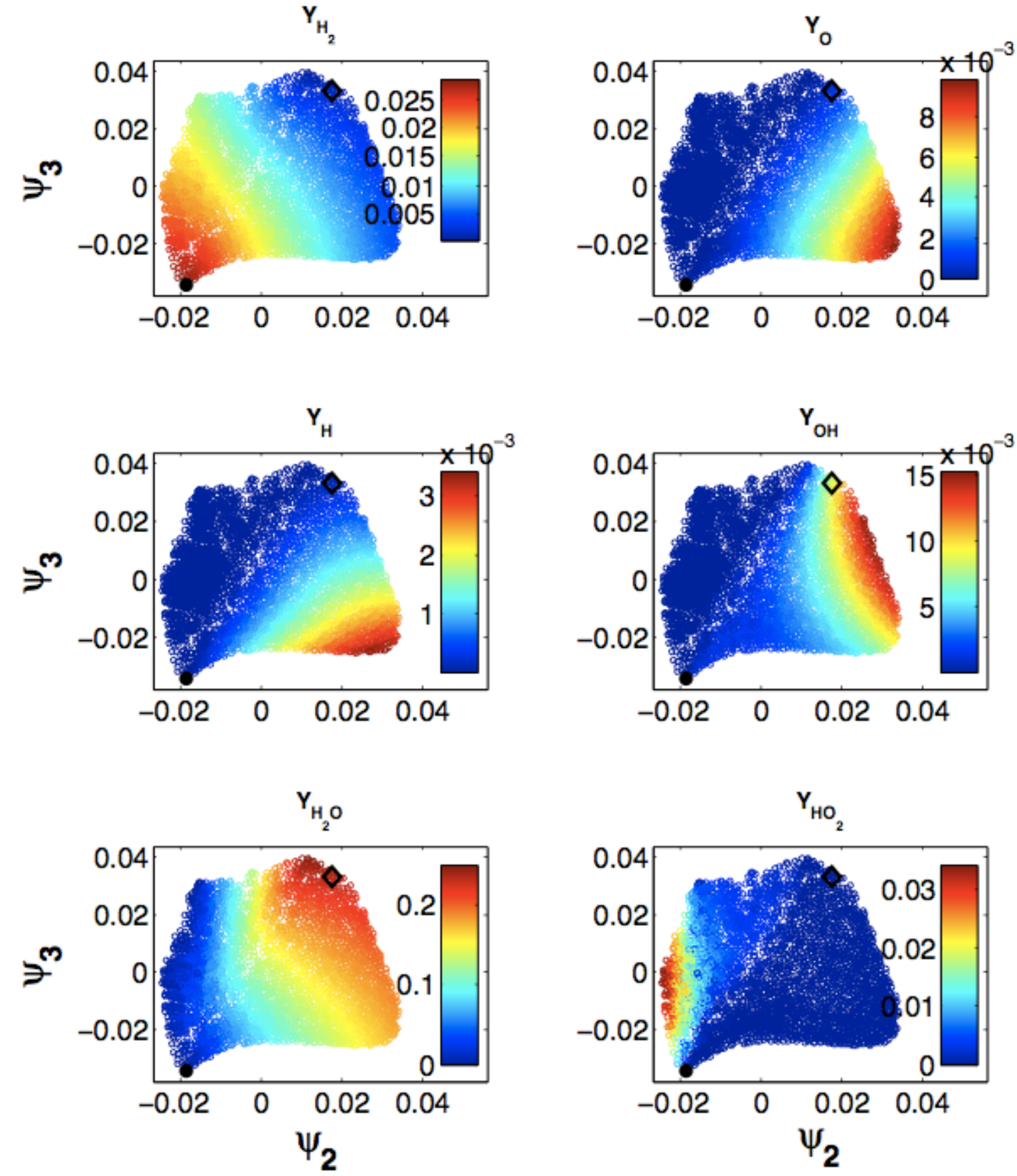}
\caption{Homogeneous reactive mixture of hydrogen and air at stoichiometric proportions with fixed enthalpy ($H=300 [kJ/kg]$) and pressure ($P=1 [bar]$). Two dimensional DMAP parameterization of 3810 points in terms of the two nontrivial leading eigenvectors $\psi_2$ and $\psi_3$ of the Markov matrix $K$. Colors represent mass fractions, while the black filled circle and the black diamond represent the fresh mixture condition and equilibrium state, respectively.} \label{figure1}
\end{figure}

\begin{figure}[t]
 \centering
 \includegraphics[scale=0.8]{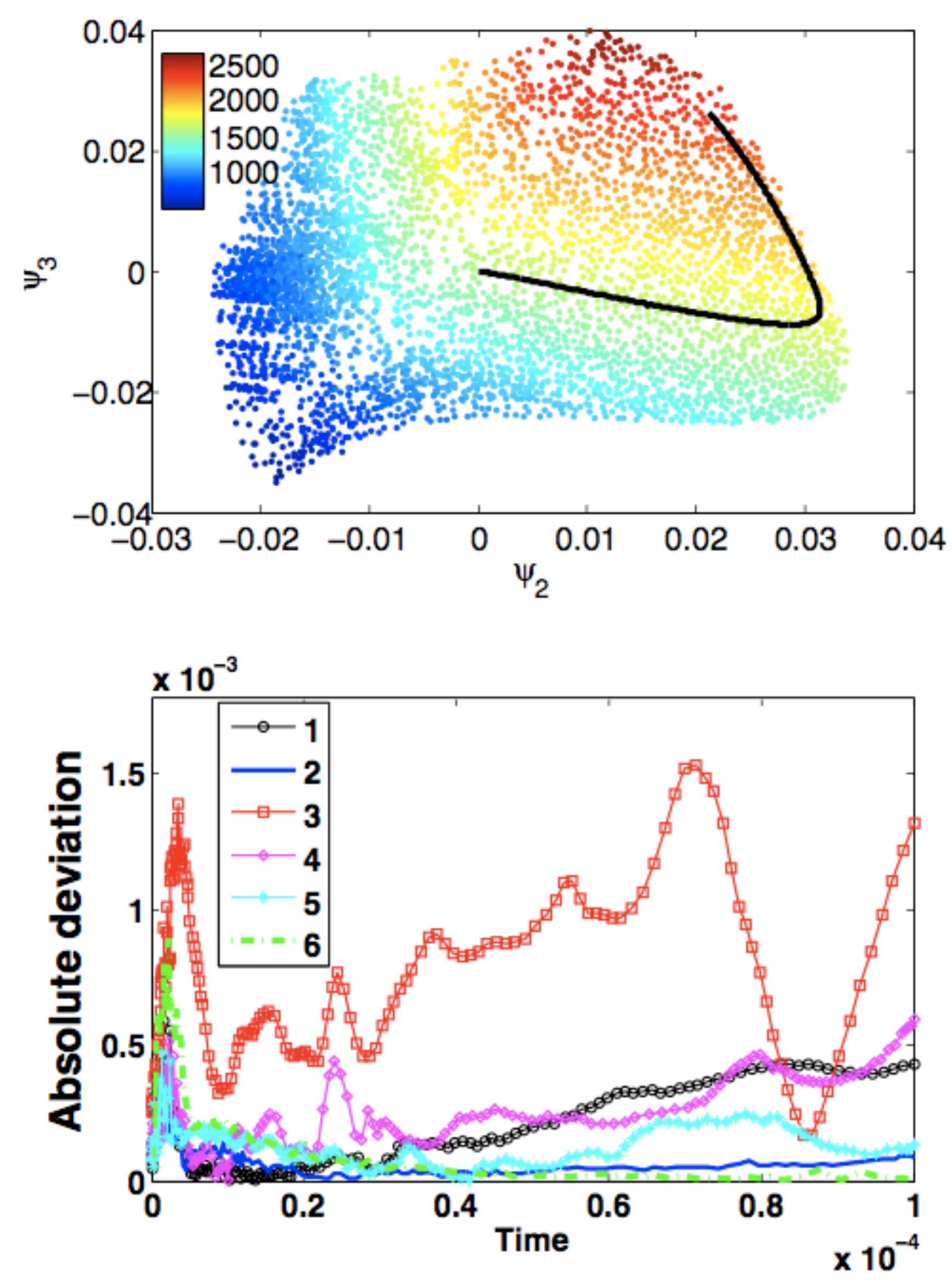}
\caption{Top: A sample detailed transient solution is shown in the plane $\psi_2-\psi_3$. Restriction is done by the Nystr\"om method, while colors refer to the temperature (Kelvin) of the gas mixture. Bottom: time evolution of the absolute deviation between detailed and reduced solution trajectories (in the reduced space) $\left\| { \psi ^{{\rm{red}}}  -  \psi ^{\det } } \right\|$. Numbers in the legend correspond to the first six methods in Table \ref{DeviationTable}.} \label{figure2}
\end{figure}

We employ our proposed approach described in Section \ref{approach} above to search for a two dimensional reduced system describing the combustion of a mixture of hydrogen and air at stoichiometric proportions under fixed total enthalpy ($H=300 [kJ/kg]$) and pressure ($P=1 [bar]$). We assume that the detailed chemical kinetics is dictated by the Li {\em et al.} mechanism \cite{Limech}, where nine chemical species ($H_2$, $N_2$, $H$, $O$, $OH$, $O_2$, $H_2O$, $HO_2$, $H_2O_2$) and three elements ($H$, $O$, $N$) are involved in the reaction. As shown in Fig. \ref{figure1}, the manifold is described by 3810 points and parameterized with respect to the two diffusion map variables $\psi_2$ and $\psi_3$. It is worth stressing that, judging from the sample density in the diffusion map space, the considered cloud of points clearly lies on a manifold with different dimensions in different regions. As expected, indeed, low temperature regions  (e.g. $T<1000$ [K]) require a larger number of reduced variables ($m >2$) to be correctly described (see Fig. \ref{figure2}) \cite{ChiavazzoKarlinPRE}. Therefore, in the example below, we only utilize {\em the portion of the manifold with high temperature (say $T>1200$ [K])}.
Coping with manifolds with varying dimension is beyond the scope of this paper, and should be addressed in forthcoming publications. We mention, however, that attempts of automatically detecting variations of the manifold dimension in the framework of diffusion maps have been also recently reported in Ref. \cite{Maggioni2011}.
\begin{figure}[t]
 \centering
 \includegraphics[scale=0.75]{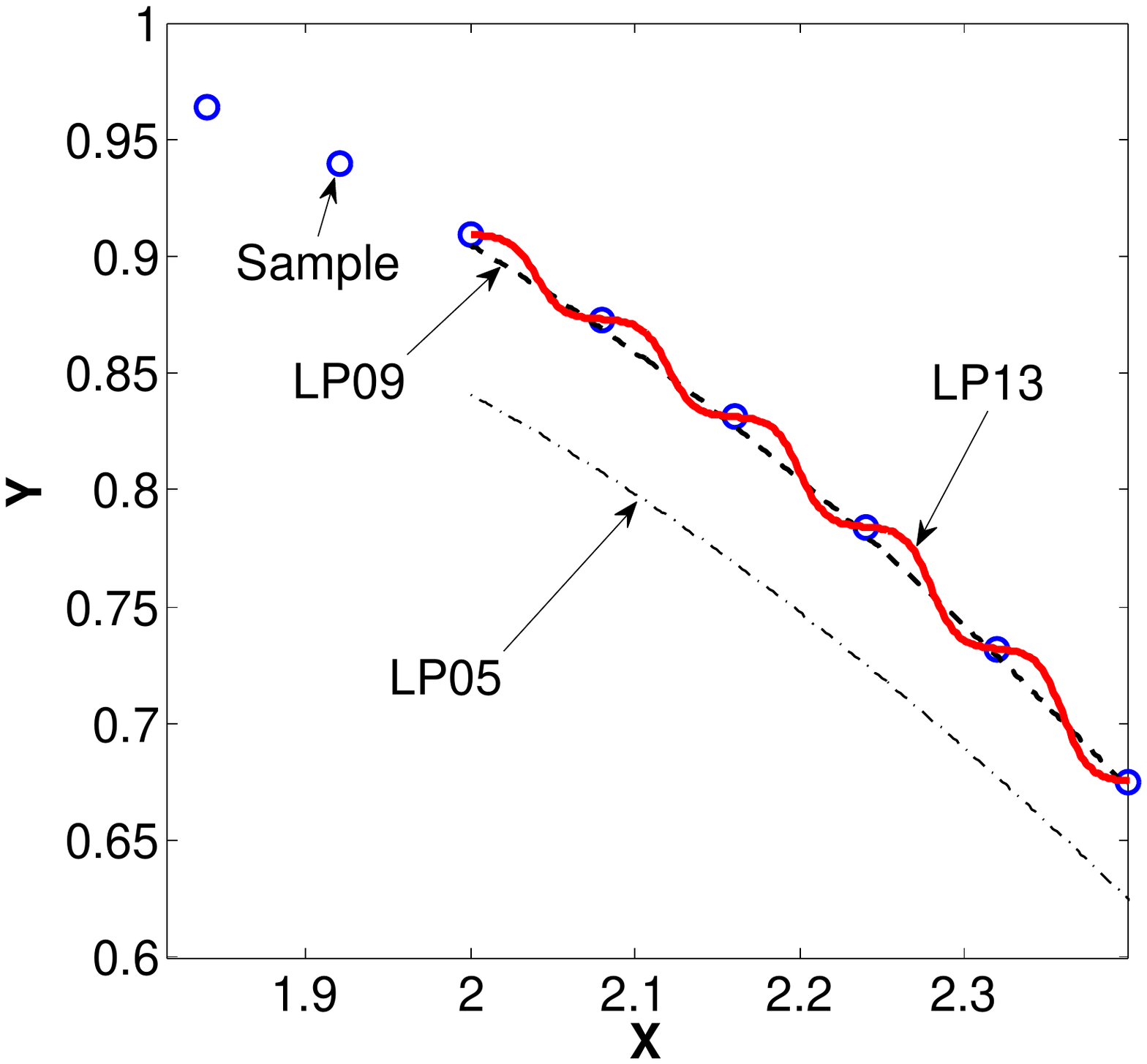}
\caption{Illustrating a possible pathology. Samples (circles) are uniformly chosen in $X$, with $Y=\sin(X)$. Laplacian Pyramids are adopted for interpolation between samples with $\sigma_0=10$. Estimated values with the finest level $l=5$, $l=9$ and $l=13$ are denoted by LP05, LP09 and LP13, respectively. At the latter level, the estimates {\em of derivatives} are no longer accurate.} \label{LPSIN}
\end{figure}
\begin{table}[htdp]
\caption{Comparison of reduced and detailed solution trajectories (with initial condition $ u=[0,0]$ and $0 \le t \le \bar t = 1 \times 10^{-4} [s]$) corresponding to several schemes implementing lifting and restriction operators (see text).  $ \left\| {\delta \psi } \right\|$ indicates the mean deviation between the reduced and detailed solution trajectory (in the reduced space): $\left\| {\delta \psi } \right\| = \bar t^{ - 1} \int_0^{\bar t} {\left\| { \psi ^{\det }  -  \psi ^{{\rm{red}}} } \right\|} dt$, with $ \psi^{\rm{det}}$, $ \psi^{\rm{red}}$  and  $\left\| { \bullet } \right\|$ denoting the restricted detailed solution, the reduced solution and the Euclidean norm, respectively. Similarly, $\left| {\delta y_\alpha} \right|$ is the mean deviation for species $\alpha$ (in the ambient space): $\left| {\delta y_\alpha } \right| = \bar t^{ - 1} \int_0^{\bar t} {\left| {y_\alpha^{\det }  - y_\alpha^{{\rm{red}}} } \right|} dt$.}
\begin{center}
{\tiny
\begin{tabular}{r|c|c|c|c|c|c|c|c|c|}
Method & $ \left\| {\delta \psi } \right\|$ & $\left| {\delta y_1 } \right|$ & $\left| {\delta y_3 } \right|$ & $\left| {\delta y_4 } \right|$ & $\left| {\delta y_5 } \right|$ & $\left| {\delta y_6 } \right|$ & $\left| {\delta y_7 } \right|$ & $\left| {\delta y_8 } \right| $ & $\left| {\delta y_9 } \right| $ \\
\hline
$1$ &   $2.28 \times 10^{-4}$ &  $2.07 \times 10^{-5}$  &  $8.52 \times 10^{-6}$ & $3.58 \times 10^{-5}$& $4.39 \times 10^{-5}$& $1.80 \times 10^{-4}$& $3.16 \times 10^{-4}$& $2.48 \times 10^{-6}$ & $1.15 \times 10^{-5}$\\
\hline
$2$ &    $5.66 \times 10^{-5}$ &  $4.09 \times 10^{-6}$ & $1.31 \times 10^{-6}$ & $5.53 \times 10^{-6}$& $1.03 \times 10^{-5}$& $3.78 \times 10^{-5}$ &$8.59 \times 10^{-5}$&  $2.18 \times 10^{-6}$ &$9.65 \times 10^{-6}$ \\
\hline
$3$  &  $8.11 \times 10^{-4}$   &  $6.90 \times 10^{-5}$   &  $2.58 \times 10^{-5}$   &  $1.04 \times 10^{-4}$   &  $1.62 \times 10^{-4}$ &   $ 6.86 \times 10^{-4}$  & $8.00 \times 10^{-4}$   &   $2.33 \times 10^{-6}$   &$ 9.98 \times 10^{-6}$ \\
\hline
$4$ &    $2.64 \times 10^{-4}$ & $2.83 \times 10^{-5}$ & $9.44 \times 10^{-6}$  &  $5.37 \times 10^{-5}$ &$1.35 \times 10^{-4}$&$2.88 \times 10^{-4}$& $2.71 \times 10^{-4}$&  $1.86 \times 10^{-6}$&$7.64 \times 10^{-6}$\\
\hline
$5$ & $1.27 \times 10^{-4}$ &  $1.26 \times 10^{-5}$ & $4.23 \times 10^{-6}$  &$2.16 \times 10^{-5}$ &$4.29 \times 10^{-5}$&$1.17 \times 10^{-4}$&$1.47 \times 10^{-4}$& $2.08 \times 10^{-6}$&$8.63 \times 10^{-6}$ \\
\hline
$6$ &     $7.31 \times 10^{-5}$ & $8.38 \times 10^{-6}$ &  $2.30 \times 10^{-6}$ & $9.15 \times 10^{-6}$&$1.78 \times 10^{-5}$&$6.76 \times 10^{-5}$& $9.46 \times 10^{-5}$& $5.70 \times 10^{-6}$ & $2.68 \times 10^{-5}$\\
\hline
$7$ &    $7.39 \times 10^{-4}$ & $7.25 \times 10^{-5}$& $2.86 \times 10^{-5}$&  $1.24 \times 10^{-4}$&$1.97 \times 10^{-4}$&$5.96 \times 10^{-4}$&  $9.94 \times 10^{-4}$&  $1.87 \times 10^{-6}$& $7.93 \times 10^{-6}$ \\
\hline
$8$  &    $8.81 \times 10^{-4}$ &  $5.90 \times 10^{-5}$ &$4.95 \times 10^{-5}$&  $1.74 \times 10^{-4}$ &$1.14 \times 10^{-4}$&$4.99 \times 10^{-4}$& $6.36 \times 10^{-4}$ & $6.27 \times 10^{-6}$ & $2.83 \times 10^{-5}$ \\
\hline
$9$  &    $0.0058$   &  $3.83 \times 10^{-4}$   &  $2.45 \times 10^{-4}$   &  $0.00107$   &  $9.66 \times 10^{-4}$ &   $ 0.0034$  & $0.0061$   &   $6.50 \times 10^{-6}$   &$ 3.75 \times 10^{-5}$ \\
\hline
$10$  &  $0.0140$   &  $0.00126$   &  $7.24 \times 10^{-4}$   &  $0.00283$   &  $0.00184$ &   $ 0.0123$  & $0.0162$   &   $2.27 \times 10^{-5}$   &$ 9.89 \times 10^{-5}$ \\
\hline
$11$  &  $8.08 \times 10^{-4}$   &  $9.03 \times 10^{-5}$   &  $4.10 \times 10^{-5}$   &  $1.63 \times 10^{-4}$   &  $1.75 \times 10^{-4}$ &   $8.33 \times 10^{-4}$  & $0.00116$   &   $2.79 \times 10^{-6}$   &$ 1.19 \times 10^{-5}$ \\
\hline
$12$  &  $0.0237$   &  $0.00331$   &  $0.00103$   &  $0.00424$   &  $0.00453$ &   $0.030$  & $0.0331$   &   $ 1.11 \times 10^{-4}$   &$ 5.45 \times 10^{-4}$ \\
\end{tabular}}
\end{center}
\label{DeviationTable}
\end{table}%

We discretized the reduced space by a $60 \times 60$ uniform Cartesian grid with $-0.025<\psi_2<0.035$ and $-0.035<\psi_3<0.04$. At every grid node, the values of the right-hand side of Eqs. (\ref{reducedODE}) (or (\ref{fastfoliationP})) are computed according to several interpolation schemes chosen form the ones described above in Section \ref{approach}, and stored in tables for later use. In particular, tables were created using the following methods:
\begin{enumerate}
\item The lifting operator consists of radial basis function interpolation with $p=3$ performed over $50$ nearest neighbors of an arbitrary point in the reduced space $u$. Restriction is done by radial basis function interpolation with $p=3$ performed over $50$ nearest neighbors of an arbitrary point in the ambient space $y$ (distances in $\Re^9$ are intended as illustrated in the Appendix). The reduced dynamical system is expressed in the form (\ref{reducedODE}).
\item The lifting operator consists of radial basis function interpolation with $p=3$ performed over $50$ nearest neighbors of an arbitrary point in the reduced space $u$. Restriction is done by the Nystr\"om method. The reduced dynamical system is expressed in the form (\ref{reducedODE}).
\item The lifting operator is based on Laplacian Pyramids up to a level $l=20$ with $\sigma_0=0.5$ over $80$ nearest neighbors of an arbitrary point in the reduced space $u$. Restriction is based on the Laplacian Pyramids up to a level $l=7$ with $\sigma_0=0.5$ over $80$ nearest neighbors of $y$. The reduced dynamical system is expressed in the form (\ref{reducedODE}). 
\item The lifting operator is based on Laplacian Pyramids up to a level $l=20$ with $\sigma_0=0.5$. Restriction is done by the Nystr\"om method. The reduced dynamical system is expressed in the form (\ref{reducedODE}).
\item  The lifting operator is based on Geometric Harmonics locally performed over $15$ nearest neighbors of an arbitrary point in the reduced space $u$. Refinements are performed until the Euclidean norm of the residual is larger than $5 \times 10^{-4}$. Restriction is done by the Nystr\"om method. The reduced dynamical system is expressed in the form (\ref{reducedODE}).
\item The lifting operator is based on Kriging performed over $8$ nearest neighbors of an arbitrary point in the reduced space $u$ (DACE package \cite{KrigSoft}, with a second order polynomial regression model, a Gaussian correlation model and parameter $\theta=10^{-3}$). Restriction is done by the Nystr\"om method. The reduced dynamical system is expressed in the form (\ref{reducedODE}). 
\item The lifting operator is based on Geometric Harmonics locally performed over $10$ nearest neighbors of an arbitrary point in the reduced space $u$. Refinements are performed until the Euclidean norm of the residual is larger than $10^{-3}$. Restriction is done using the Nystr\"om method. The reduced dynamical system is expressed in the form (\ref{reducedODE}).
\item The lifting operator is based on Kriging performed over $8$ nearest neighbors of an arbitrary point in the reduced space $u$ (DACE package \cite{KrigSoft}, with second order polynomial regression model, a Gaussian correlation model and parameter $\theta=10^{-3}$). Restriction is done by the Nystr\"om method. The reduced dynamical system is expressed in the form (\ref{fastfoliationP}).
\item The lifting operator is based on Kriging performed globally over all samples (package \cite{KrigSoft}, with a second order polynomial regression model, a Gaussian correlation model and parameter $\theta=13$). Restriction is done by the Nystr\"om method. The reduced dynamical system is expressed in the form (\ref{fastfoliationP}).
\item The lifting operator is based on the Laplacian Pyramids up to a level $l=20$ with $\sigma_0=0.5$ over $80$ nearest neighbors of an arbitrary point in the reduced space $u$. Restriction is based on the Laplacian Pyramids up to a level $l=3$ with $\sigma_0=0.5$ over $80$ nearest neighbors of an arbitrary point in the ambient space $y$. The reduced dynamical system is expressed in the form (\ref{reducedODE}). 
\item The lifting operator is based on the Laplacian Pyramids up to a level $l=20$ with $\sigma_0=0.5$ over $80$ nearest neighbors of an arbitrary point in the reduced space $u$. Restriction is based on the Laplacian Pyramids up to a level $l=9$ with $\sigma_0=0.5$ over $80$ nearest neighbors of an arbitrary point in the ambient space $ y$. The reduced dynamical system is expressed in the form (\ref{reducedODE}).
\item The lifting operator is based on the Laplacian Pyramids up to a level $l=20$ with $\sigma_0=0.5$ over $80$ nearest neighbors of an arbitrary point in the reduced space $u$. Restriction is based on Laplacian Pyramids up to a level $l=12$ with $\sigma_0=0.5$ over $80$ nearest neighbors of an arbitrary point in the ambient space $y$. The reduced dynamical system is expressed in the form (\ref{reducedODE}).
\end{enumerate}
%
%
Each of the above tables was utilized for providing the systems (\ref{reducedODE}) and (\ref{fastfoliationP}) with a closure, where rates of reduced variables are efficiently retrieved via bi-variate interpolation in diffusion map space. In Fig. \ref{figure2} a sample trajectory (starting from $u=[0,0]$) is reported in the top part, while the Euclidean norm of the absolute deviation between the reduced and detailed solution (in the $\psi_2 - \psi_3$ plane) is reported in the lower part of the figure as a function of time. A more detailed comparison is reported in the Table \ref{DeviationTable}. In our (not optimized) implementation, all trajectories are computed by the Matlab's solver {\bf ode45}, with the reduced system showing a speedup of roughly four times compared to the detailed one.

In terms of accuracy, we found that the best performances are achieved combining a local lifting operator (e.g. interpolation/extension over nearest neighbors) with the Nystr\"om method for restriction. For instance, we notice that a proper combination between radial basis function interpolation (for lifting) and Nystr\"om extension may offer excellent accuracy (in terms of deviation errors $\left\| {\delta \psi } \right\|$ and $\left| {\delta y_i} \right|$), as showed in Table \ref{DeviationTable} for the solution trajectory  in Fig. \ref{figure2}. Clearly, radial basis functions are simpler to implement and require less computational resources compared to other approaches such as Kriging and Geometric Harmonics. We should stress, though, that the latter techniques present similar performances and are certainly to be preferred in cases where (unlike Figs. \ref{figure1} and \ref{figure2}) samples are not uniformly distributed (i.e. sample clustering). Moreover, we observe that approaches based on Laplacian Pyramids (for restriction) present poorer performances even with large values of $l$.
An explanation for this is a possible inaccurate estimate of the derivatives at the right-hand side of the reduced dynamical system, which we attempt to illustrate
through the caricature in Fig. \ref{LPSIN}.
We finally find that solutions to the system (\ref{fastfoliationP}) typically lead to larger errors compared to those obtained solving (\ref{reducedODE}).

\begin{figure}[t]
 \centering
 \includegraphics[scale=0.65]{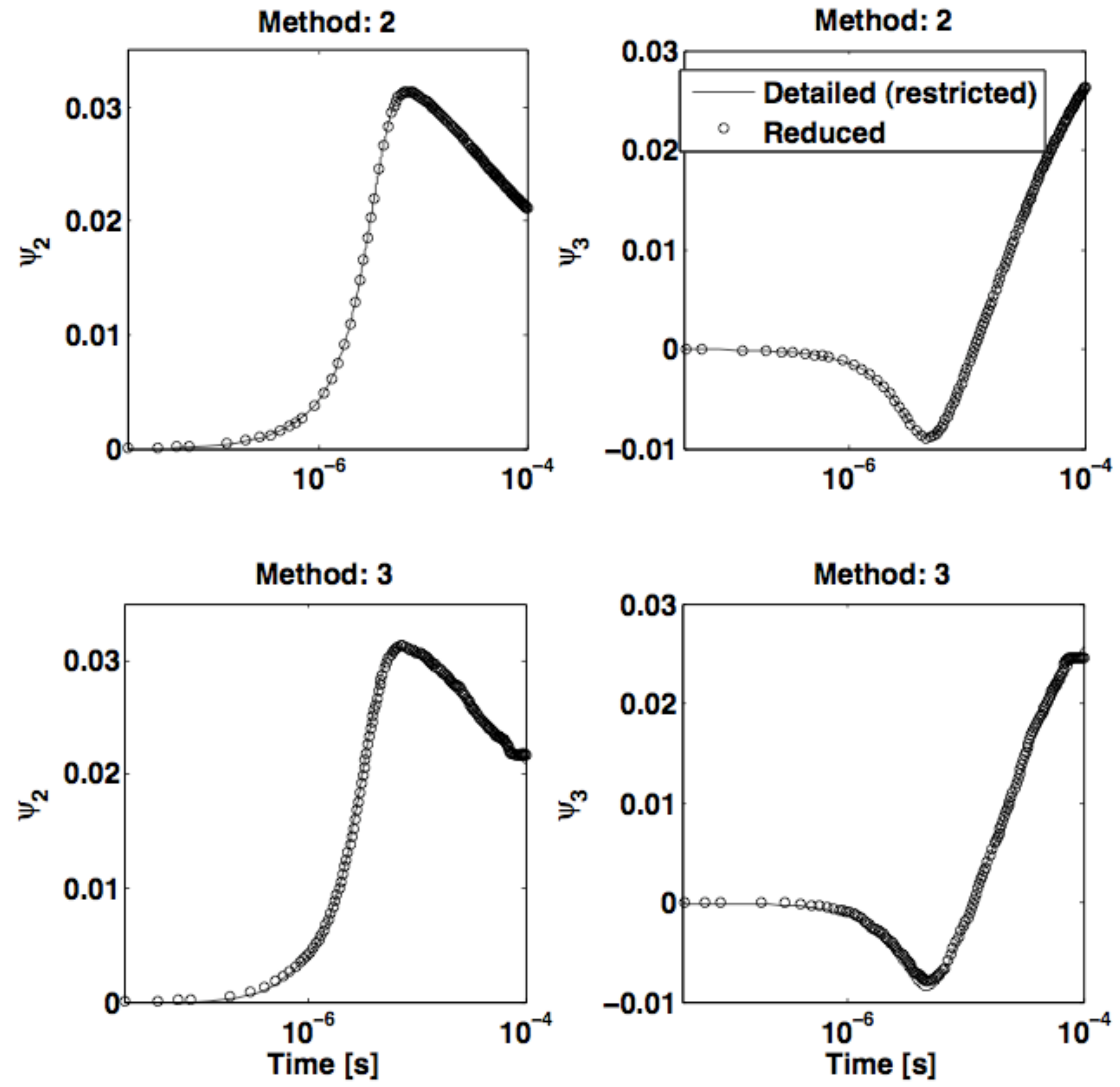}
\caption{Time evolution of the two diffusion map variables along the solution trajectory of Fig. \ref{figure2} as obtained by Method 2 (top) and by Method 3 (bottom) (see Table 1). The initial condition in the diffusion maps space $[0,0]$ is first lifted into $\Re^{9}$ and then relaxed towards the equilibrium point by the detailed kinetics (\ref{detailedODE}) using the readily available Matlab solver {\bf ode45}. The latter time series is afterwards restricted to the diffusion maps space and reported with a continuous line. Symbols denote the corresponding solution directly obtained in the reduced space by solving the system (\ref{reducedODE}) by the same Matlab solver {\bf ode45}.} \label{TSdmaps}
\end{figure}

\begin{figure}[t]
 \centering
 \includegraphics[scale=0.75]{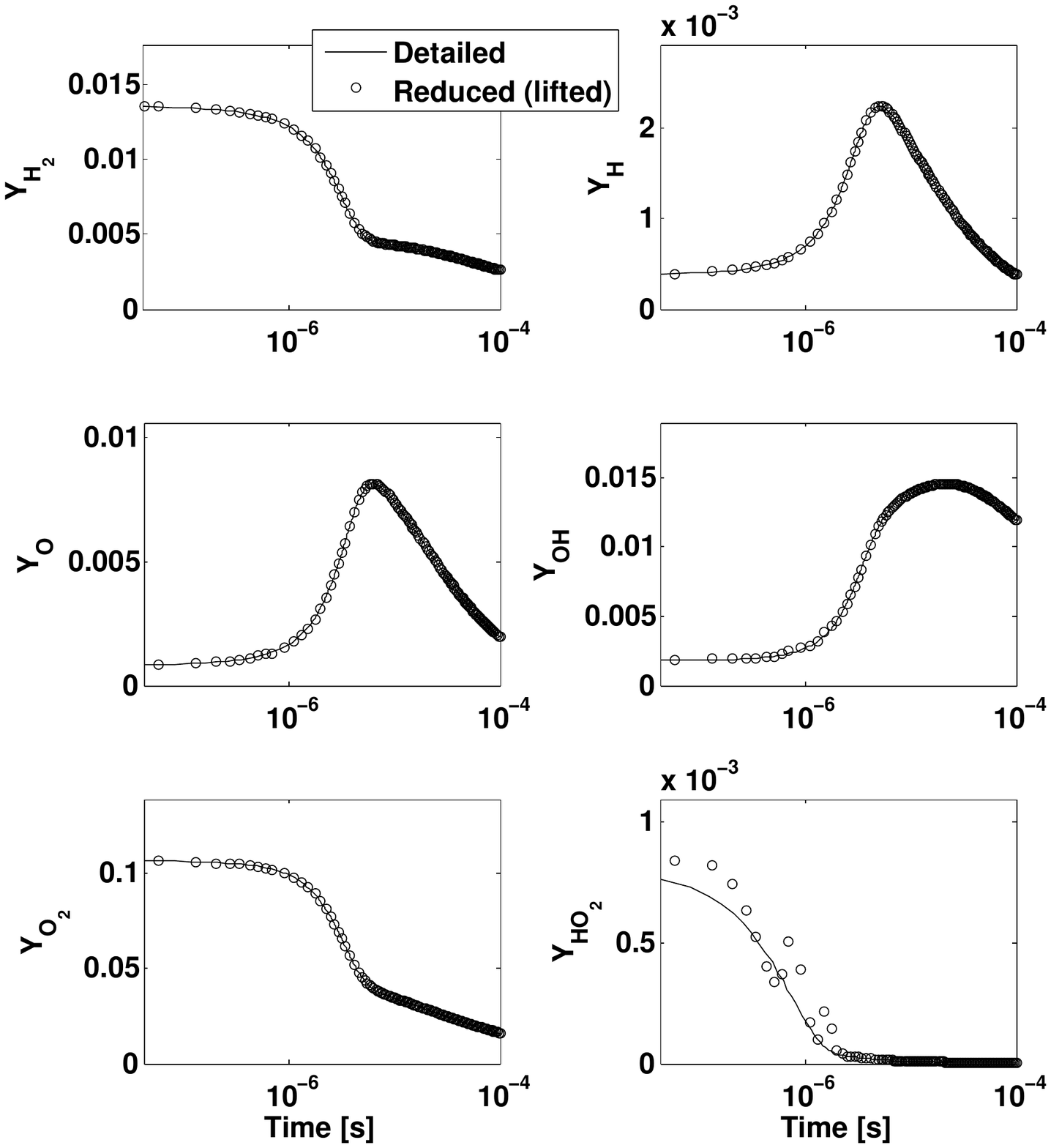}
\caption{The initial condition in the diffusion maps space $[0,0]$ is first lifted into $\Re^{9}$ and then relaxed towards the equilibrium point by the detailed kinetics (\ref{detailedODE}) using the readily available Matlab solver {\bf ode45} (continuous line). Symbols report the corresponding time series as obtained by lifting the reduced solution at the top of Fig. \ref{TSdmaps} (i.e. method 2).} \label{TSspeciesM02}
\end{figure}

\begin{figure}[t]
 \centering
 \includegraphics[scale=0.75]{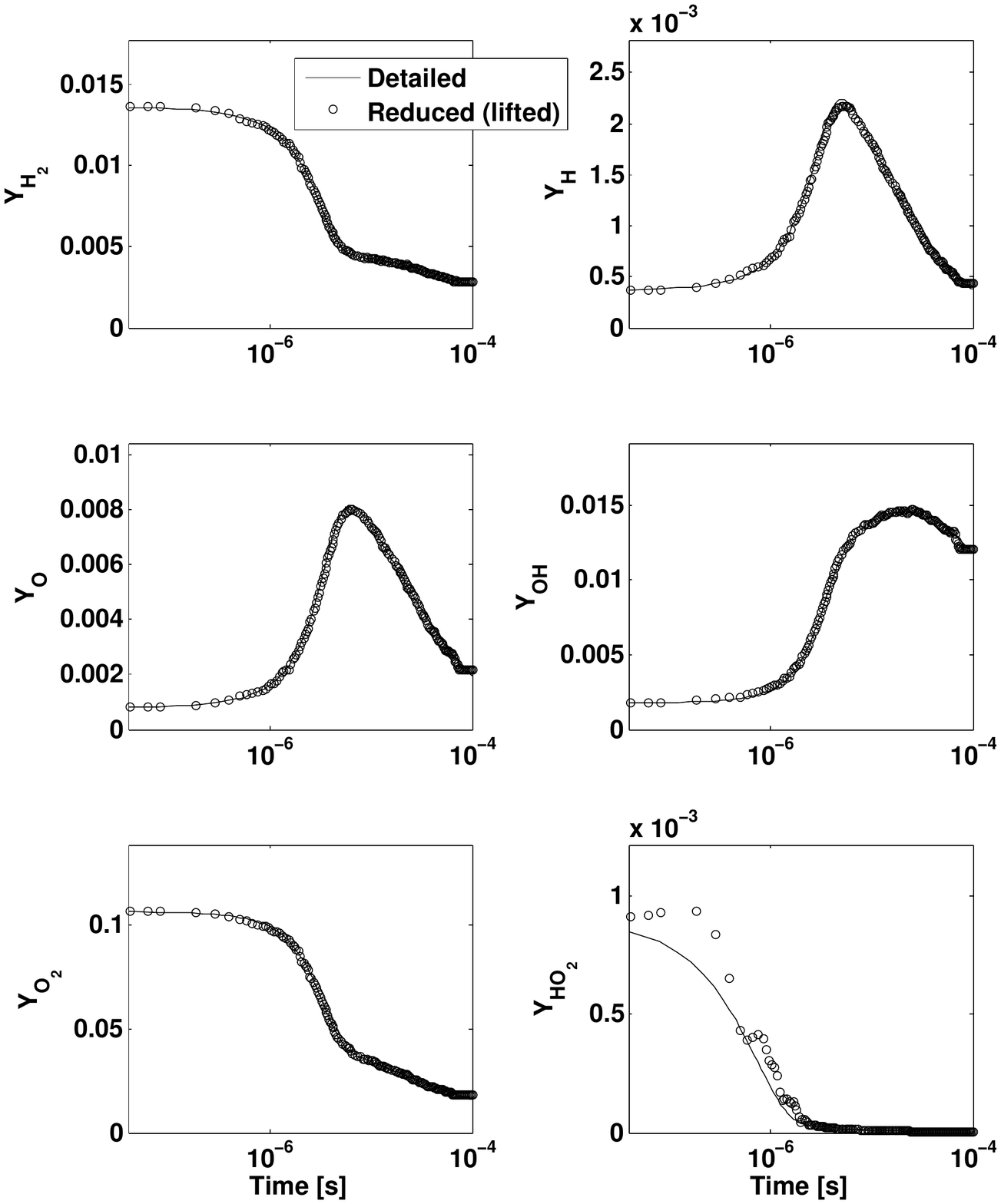}
\caption{The initial condition in the diffusion maps space $[0,0]$ is first lifted into $\Re^{9}$ and then relaxed towards the equilibrium point by the detailed kinetics (\ref{detailedODE}) using the readily available Matlab solver {\bf ode45} (continuous line). ). Symbols report the corresponding time series as obtained by lifting the reduced solution at the bottom of Fig. \ref{TSdmaps} (i.e. method 3).} \label{TSspeciesM03}
\end{figure}

%
%
%
For completeness, in Fig. \ref{TSdmaps} we report the time series of the diffusion maps variables as obtained by the methods 2 and 3 in the Table \ref{DeviationTable}, as well as the restriction of the corresponding detailed solution. Moreover, in Figs. \ref{TSspeciesM02} and \ref{TSspeciesM03} a comparison of the time series in the detailed space is reported as obtained by reconstruction of the states in $\Re^9$ from the reduced solutions in Fig. \ref{TSdmaps}.
%
\section{Conclusions}\label{conclusion}
In this work, we showed that the diffusion maps (DMAP) technique is a promising tool for extracting a global parameterization of low-dimensional manifolds arising in combustion problems. 
Based on the slow variables automatically  identified by the process, a reduced dynamical system can be obtained and solved. 
Both lifting and restriction operators (i.e. mapping of any point in the region of interest of the reduced space into the full space and vice-versa) lie at the heart of such an approach. 
To construct these operators, methods for extending empirical functions only known at scattered locations must be employed, and we have tested several.

For chemical kinetics governing a non-isothermal reactive gas mixture of hydrogen and air, a comparison is carried out on the basis of the deviation error between sample detailed solutions and the corresponding reduced ones in both the full and reduced spaces. 
Several combinations of interpolation schemes were implemented in the procedure restrictions/liftings, with the reduced rates ${{d u} \mathord{\left/{\vphantom {{d u} {dt}}} \right.\kern-\nulldelimiterspace} {dt}}$ pre-computed and stored in tables to be utilized at a later time for providing the system (\ref{reducedODE}) with a closure.
In the considered case, approaches based on a local lifting operator (i.e. interpolation/extension over nearest neighbors) combined with the Nystr\"om method (for restriction) have shown superior performances in terms of accuracy in recovering the (longer-time) transient dynamics of the
detailed model.

While the feasibility of the presented approach has been demonstrated here, a number of open issues remain. In particular, future studies should focus on computationally efficient implementations of the method without pre-tabulation, since handling tables at high dimensions (say $m>4$) becomes computationally complex. Moreover, as demonstrated also in the presented combustion example, the method should be able to cope with manifolds whose dimension possibly varies across distinct regions of the phase-space; how to consistently express and solve reduced systems across manifolds with disparate dimensions remains out of reach for the present method,  requiring further investigation.

\section{Appendix}\label{appendix}
Due to a disparity of the magnitudes of species concentrations, $d_{ij}$ is taken as the Euclidean distance between properly rescaled points $ y'_i$ and $y'_j$, with $y'_i = R y_i$ using the fixed diagonal matrix $R=\{ r_{\beta \beta}\}$, $r_{\beta \beta}=1/max(y_\beta)$. Here, $max(y_\beta)$ represents the largest $\beta$-th coordinate among all available samples.

\acknowledgements{Acknowledgements}
E.C. acknowledges partial support of the US-Italy Fulbright Commission and the Italian Ministry of Research (FIRB grant RBFR10VZUG). I.G.K. and C.W.G. gratefully acknowledge partial support by the US DOE. C.J.D. acknowledges support by the US Department of Energy Computation Science Graduate Fellowship (grant number DE-FG02-97ER25308).

\bibliographystyle{mdpi}
\makeatletter
\renewcommand\@biblabel[1]{#1. }
\makeatother

%
%
\end{document}